\documentclass[reqno,12pt, oneside]{amsart}

\usepackage{enumerate}
\usepackage[nobysame]{amsrefs}
\usepackage{enumitem}
\usepackage{amssymb}
\usepackage{amsmath}
\usepackage{amsthm}
\usepackage{mathtools}
\usepackage{amscd}
\usepackage{microtype}
\usepackage[right=2cm,left=2cm, top=2.5cm, bottom=2.5cm,includehead]{geometry}
\usepackage{amsfonts}
\usepackage{booktabs}
\usepackage{hyperref}

\usepackage{tgpagella}
\usepackage{eulervm}

\usepackage{chngcntr}
\counterwithin{table}{section}

\usepackage[backgroundcolor=lightgray]{todonotes}

\theoremstyle{plain}
\newtheorem{theorem}{Theorem}
\newtheorem{corollary}[theorem]{Corollary}
\newtheorem{lemma}[theorem]{Lemma}
\newtheorem{proposition}[theorem]{Proposition}
\newtheorem{theorema}{Theorem}

\theoremstyle{remark}
\newtheorem{remark}[theorem]{Remark}
\newtheorem{example}[theorem]{Example}

\theoremstyle{definition}

\DeclarePairedDelimiter{\abs}{\lvert}{\rvert}

\DeclarePairedDelimiter{\norm}{\lVert}{\rVert}

\DeclarePairedDelimiterX{\dset}[2]{\lbrace}{\rbrace}{#1\;\delimsize|\;#2}

\DeclareMathOperator{\Id}{Id}

\DeclareMathOperator{\var}{var}

\newcommand{\ball}{{\overline{B}}}

\DeclareMathOperator{\Lip}{Lip}

\title[Composition operators in $bv_p$-spaces, part I]{Composition operators in $bv_p$-spaces,\\ part I: acting conditions and boundedness}

\author{Daria Bugajewska}
\address[D.~Bugajewska]{Faculty of Mathematics and Computer Science\\
  Adam Mickiewicz University in Pozna\'n\\
  ul.\ Uniwersytetu Pozna\'nskiego 4\\
  61-614 Pozna\'n\\
  Poland}
\email[D.~Bugajewska]{dbw@amu.edu.pl}

\author{Piotr Kasprzak}
\address[P. Kasprzak]{Department of Nonlinear Analysis and Applied Topology\\
Faculty of Mathematics and Computer Science\\
  Adam Mickiewicz University in Pozna\'n\\
  ul.\ Uniwersytetu Pozna\'nskiego 4\\
  61-614 Pozna\'n\\
  Poland}
\email[P.~Kasprzak]{kasp@amu.edu.pl}


\date{\today}

\begin{document}
\begin{abstract}
The aim of this paper is to give the answer to the problem of characterization of acting conditions (necessary as well as sufficient) for composition operators in some sequence spaces. We also characterize their boundedness and local boundedness. We focus on composition operators acting to or from the space $bv_p(E)$ of all sequences of $p$-bounded variation; here $p\geq 1$ and $E$ is a normed space.
\end{abstract}

\maketitle

\section{Introduction}

Superposition operators (also known as Nemytskii operators) play an important role in numerous areas of mathematical and nonlinear analysis. This is because many problems involving integral or differential equations can be solved using the fixed-point formulation $x=F(x)$, where $F=K\circ S_f$ with $K$ often being an integral operator generated by a Green's function and $S_f$ being a superposition operator.
 
The theory of superposition operators in various function or sequence spaces is both vast and rich. However, the most significant results are scattered throughout the literature. To the best of our knowledge, there is only one comprehensive source in English dedicated entirely to such nonlinear mappings, and that is the celebrated monograph by Appell and Zabrejko~\cite{AZ}. 

In the context of classical sequence spaces, a notable paper is the work of Dedagich and Zabrejko (see~\cite{DZ}). This article discusses certain function-theoretic and topological properties of superposition operators, primarily focusing on $l^p$-spaces. Other works worth mentioning are \cites{BAM, K, LF, Pl, Ro, Sh}. 
Even more surprising is the fact that the study of superposition operators in another classical space (or, ``special space'' as Dunford and Schwarz call it in their book -- see~\cite{DS}*{p.~239}), that is the space $bv_p(E)$ of sequences of bounded $p$-variation, has only begun very recently. One of the first papers addressing this topic was the work by Karami et al.~\cite{KFA}. As a side note, let us mention that the parallel theory of superposition operators in spaces of functions of bounded variation of various types is well-developed (for more details see, for example, monographs~\cite{ABM}*{Chapters~5 and~6}, \cite{DN}*{Chapter~6} and~\cite{reinwand}*{Chapter~5}, or articles\cites{BBKM, BCGS, AGV, AMS, MM, M}.  

There are at least two reasons why the spaces $bv_p(E)$ are especially interesting. Although they are linearly isometric to $E\times l^p(E)$, they contain unbounded sequences if $p>1$. Moreover, it can be shown that $bv_1(\mathbb R)$ can be naturally interpreted as the dual to the space of all non-absolutely convergent real series (see~\cite{DS}*{Exercise~12, p.~339}).

In this paper we focus on composition operators, that is, those superposition operators which are generated by functions independent of the ``time'' variable; for precise definitions see Section~\ref{sec:composition_operators}. The main goal is to answer the fundamental question of the theory and to describe necessary and sufficient conditions under which the given composition operator acts between $bv_p(E)$ and one of the spaces $c_0(E)$, $c(E)$, $l^q(E)$ or $bv_r(E)$; here $p,r\in [1,+\infty)$ and $q \in [1,+\infty]$. We also completely characterize (local) boundedness of such operators.

The motivation for our study is the aforementioned paper~\cite{KFA}. It seems that the results established in that article may have some flaws, to say the least. The main concern is the fact that in~\cite{KFA} the properties of superposition operators are actually expressed in terms of the operators themselves. Not only is it not in line with a common practice in the literature of the subject, but it also makes it nearly impossible to apply those results. Additionally, there are some doubts regarding the validity of the presented proofs. (For more detailed discussion of the issues concerning~\cite{KFA} please refer to Section~\ref{sec:Karami}.)  

Let us emphasize that, in contrast to~\cite{KFA}, we describe the acting conditions and (local) boundedness of composition operators to and from $bv_p$-spaces in terms of their generators. Furthermore, we employ completely different methods than those used in~\cite{KFA}. Finally, rather than focusing on the real-valued case, we consider sequences with values in arbitrary normed spaces. As shown in the recent article~\cite{BK}, the transition from the finite to infinite-dimensional setting significantly impacts the properties of composition/superposition operators.    

The paper is organized as follows. In Section 2, not only do we introduce notation and recall basic facts concerning classical sequence spaces, but we also study several classes of H\"older continuous mappings, discussing their properties and relationships. These maps will play a key role in the subsequent parts of the article. Section 3 is dedicated to characterizing acting conditions for composition operators. We focus on those which act between $bv_p(E)$ and one of the spaces $c_0(E)$, $c(E)$, $l^q(E)$ or $bv_r(E)$, where $p,r\in [1,+\infty)$ and $q \in [1,+\infty]$. Section 4 deals with the (local) boundedness of the composition operators considered in the previous parts of the paper. Finally, in the last fifth section, we provide a brief comparison of properties of composition operators acting in $bv_p$- and $BV_p$-spaces. We also discuss the results obtained in~\cite{DZ} and~\cite{KFA} and compare them with our findings.

\section{Preliminaries}

The aim of this section is twofold. First, we will introduce notation and recall some basic facts concerning classical sequence spaces and composition operators. And then, we will discuss various classes of generators. 

\subsection{Notation}
By $\mathbb N$ we will denote the set of all positive integers. Let $E$ be a normed space over the field of either real or complex numbers.
The open and closed balls in $E$ with center $u \in E$ and radius $r > 0$ will be denoted by $B_E(u,r)$ and $\ball_E(u, r)$, respectively. Given two sets $A$ and $B$, we will write $A\subseteq B$ if $A$ is included in $B$, and $A \subset B$ if the inclusion is strict. Furthermore, for a non-empty set $A$ we will denote its characteristic function by $\chi_A$. Also, if $f \colon E \to E$ is a map and $A\subseteq E$, then $f|_A$ will stand for the restriction of $f$ to the set $A$. Finally, we will say that $f$ is bounded on $A$, if $\sup_{v \in A}\norm{f(v)}<+\infty$. Regarding summation, we will assume that the sum over the empty set or when the upper summation index exceeds the lower one always gives zero. Lastly, it is worth noting that we will use the same symbol ``$0$'' to denote the zero scalar, the zero vector (the origin of a normed space) and the zero sequence. This will lead to no confusion and will keep down the amount of symbols to be used throughout the paper. 

\subsection{Sequence spaces}
In this note, we will adopt functional notation to represent elements of sequence spaces, that is, $x(n)$ will stand for the $n$-th term of the sequence $x$. On the other hand, when talking about sequences of scalars, elements of normed spaces, or even other sequences, we will use the standard notation $(u_n)_{n \in \mathbb N}$.

Let $p \in [1,+\infty)$ and let $(E,\norm{\cdot})$ be a normed space over the field of either real or complex numbers. (Note that we do not assume here that $E$ is complete.) By $l^p(E)$ we will denote the space of all $E$-valued sequences that are absolutely summable with $p$-th power, endowed with the norm
\[
 \norm{x}_{l^p} := \Biggl(\sum_{n=1}^\infty \norm{x(n)}^p \Biggr)^{\frac{1}{p}}.
\]  
Furthermore, by $l^\infty(E)$ we will denote the normed space of all $E$-valued bounded sequences with the supremum norm $\norm{x}_\infty := \sup_{n \in \mathbb N}\norm{x(n)}$, and by $c(E)$ and $c_0(E)$ its normed subspaces of all, respectively, convergent and null sequences. In several counterexamples we will also use yet another subspace of $l^{\infty}(E)$, that is $c_{00}(E)$, which consists of all $E$-valued sequences with only finitely many non-zero terms. The main object of our study will be the space $bv_p(E)$ of all sequences $x \colon \mathbb N \to E$ such that $\sum_{n=1}^\infty \norm{x(n+1)-x(n)}^p<+\infty$. The default norm on $bv_p(E)$ is 
\[
 \norm{x}_{bv_p}:=\norm{x(1)}+ \Biggl(\sum_{n=1}^\infty \norm{x(n+1)-x(n)}^p\Biggr)^{\frac{1}{p}}.
\]
If $E$ is a Banach space, all the sequence spaces considered above except $c_{00}(E)$ are complete.

At this point lest us also recall some inclusions holding between the aforementioned sequence spaces. For $1 \leq p < q <+\infty$ we have
\[
 l^p(E) \subset l^q(E) \subset c_0(E) \subset c(E) \subset l^\infty(E) \quad \text{and} \quad l^p(E) \subset bv_p(E) \subset bv_q(E). 
\]
Further, $bv_1(E) \subset l^{\infty}(E)$ and $\norm{x}_\infty \leq \norm{x}_{bv_1}$ for $x \in bv_1(E)$. But if $p>1$, then neither $bv_p(E)$ is included in $l^{\infty}(E)$, nor $l^{\infty}(E)$ in $bv_p(E)$. (To see this it suffices to take $E:=\mathbb R$ and consider the sequences $(1,1+\frac{1}{2},1+\frac{1}{2}+\frac{1}{3},\ldots)$ and $(0,1,0,1,0,1,\ldots)$, respectively.) In particular, $bv_p(E) \not\subseteq c(E)$ for $p>1$. Finally, if the normed space $E$ is complete, then $bv_1(E)\subset c(E)$ . What may come as a surprise is the fact that completeness of $E$ is essential for the above inclusion to hold. To see this let us take the sequence $x \colon \mathbb N \to c_{00}(\mathbb R)$ defined by the formula $x(n):=(1,\frac{1}{2^2},\frac{1}{3^2},\ldots,\frac{1}{n^2},0,0,\ldots)$. Then, $\norm{x(n+1)-x(n)}_{\infty} = (n+1)^{-2}$ for $n \in \mathbb N$, which implies that $x \in bv_1(c_{00}(\mathbb R))$. However, $x \notin c(c_{00}(\mathbb R))$. (We refer the reader to~\cite{BA} for more details and some proofs of the facts related to the above inclusions for $bv_p$-spaces in the case $E:=\mathbb R$.) 

In sequence spaces one can define very natural projections. Given $n \in \mathbb N$ by $P_n$ we will denote the projection onto the first $n$ components. Also, we will set $Q_n:=\Id-P_n$; here $\Id$ stands for the identity map. In other words, for a sequence $x \colon \mathbb N \to E$ we have $P_n(x)=(x(1),x(2),\ldots,x(n),0,0,\ldots)$ and $Q_n(x)=(0,\ldots,0,x(n+1),x(n+2),\ldots)$.

\subsection{Space of functions of Wiener bounded variation}

In the last section we are going to compare the results obtained for $bv_p$-spaces with those concerning the spaces of functions of Wiener bounded $p$-variation that are already known in the literature. Therefore, let us briefly recall the definition of the space $BV_p([a,b], E)$. As before, let $p$ be a fixed real number in the interval $[1,+\infty)$ and let $E$ be a normed space. The (possibly infinite) quantity
\[
 \var_p x:=\sup \sum_{i=1}^n\norm{x(t_{i})-x(t_{i-1})}^p,
\]
 where the supremum is taken over all finite partitions $a=t_0<t_1<\cdots<t_n=b$ of $[a,b]$, is called the \emph{$p$-variation} (or, \emph{Wiener variation}) of the function $x$ over the interval $[a,b]$. If $\var_p x<+\infty$, then we say that $x$ is of \emph{bounded $p$-variation}. The linear space of all such maps, denoted by $BV_p([a,b],E)$, is a normed space when endowed with the norm $\norm{x}_{BV_p}:=\norm{x(a)} + (\var_p x)^{1/p}$. (More information on the Wiener variation, or its generalization -- the Young variation, can be found in~\cite{DN}. The real-valued case is also thoroughly discussed in~\cite{ABM}.)

\subsection{Composition operators}
\label{sec:composition_operators}

Let $E$ be a normed space and let $\Omega$ be a non-empty set. Moreover, let $f \colon \Omega \times E \to E$ be a map. The \emph{superposition operator} $S_f$,  corresponding to $f$, to a function $x \colon \Omega \to E$ (which usually belongs to some given set or space) assigns the mapping $S_f(x)$ defined by $S_f(x)(t)=f(t,x(t))$ for $t \in \Omega$. In the literature, superposition operators are often also referred to as \emph{Nemytskii operators}. The function $f$ is called the \emph{generator} of $S_f$. 

If $f$ does not depend on the first ``time'' variable, that is $f \colon E \to E$, we talk about \emph{autonomous superposition} or \emph{Nemytskii operators}. In this paper, just for simplicity, we will use their yet another name -- \emph{composition operators}. Also, instead of $S_f$ we will denote them by $C_f$. As we will mainly deal with sequence spaces, throughout the paper (except for the last section) we will assume that $\Omega:=\mathbb N$. In Section~\ref{sec:comparison} we will consider the case $\Omega:=[a,b]$, too.

For more details on composition and superposition operators we refer the reader to the celebrated monograph by Appell and Zabrejko~\cite{AZ}.

\subsection{Generators}
\label{sec:generators}

Generators of composition operators in $bv_p$-spaces with nice (topological) properties revolve around H\"older continuous mappings. Because in the literature there is no consensus on how to name various classes of such maps, we decided to recall all the necessary definitions. Moreover, we will briefly discuss their properties and relationships. Note, however, that we will not go very deep into details, as H\"older maps are not the main object of our study.

Let $E$ be a normed space and let $\alpha\in (0,1]$. A map $f\colon E \to E$ is called
\begin{enumerate}[label=\textup{(\alph*)}]
 \item \emph{H\"older continuous with exponent $\alpha$}, if there exists $L\geq 0$ such that $\norm{f(u)-f(w)}\leq L\norm{u-w}^\alpha$ for any $u,w \in E$; the set of all such maps will be denoted by $\Lip^\alpha(E)$,

 \item \emph{locally H\"older continuous in the stronger sense with exponent $\alpha$}, if there exist $\delta>0$ and $L\geq 0$ such that $\norm{f(u)-f(w)}\leq L\norm{u-w}^\alpha$ for all $u,w \in E$ with $\norm{u-w}\leq \delta$; the set of all such maps will be denoted by $\Lip^\alpha_{\text{strg}}(E)$,

 \item \emph{H\"older continuous on bounded sets with exponent $\alpha$}, if for every $r>0$ there exists $L_r \geq 0$ such that $\norm{f(u)-f(w)}\leq L_r\norm{u-w}^\alpha$ for all $u,w \in \ball_E(0,r)$; the set of all such maps will be denoted by $\Lip^\alpha_{\text{bnd}}(E)$,

 \item \emph{H\"older continuous on compact sets with exponent $\alpha$}, if for every compact subset $K$ of $E$ there exists $L_K \geq 0$ such that $\norm{f(u)-f(w)}\leq L_K\norm{u-w}^\alpha$ for all $u,w \in K$; the set of all such maps will be denoted by $\Lip^\alpha_{\text{comp}}(E)$.
\end{enumerate}
Very often, if the exponent $\alpha$ is known from the context and there is no chance of confusion, we will omit the part ``with exponent'' and write, for example, ``$f$ is H\"older continuous''  instead of ``$f$ is H\"older continuous with exponent $\alpha$.'' Furthermore, to maps which are H\"older continuous (in any sense) with exponent $\alpha=1$ we will refer simply as Lipschitz continuous.

\begin{remark}\label{rem:holder_alpha_greater_than_1}
Of course, the above definitions are well-posed for any $\alpha >0$. However, it can be easily shown that if $\alpha>1$ all the above classes of functions consists of constant maps only.
\end{remark}

\begin{remark}
Most of the above notions are classical. The motivation to consider mappings satisfying local H\"older continuity in the stronger sense comes from the paper~\cite{KFA}, where a similar condition was used when dealing with superposition operators acting between different $bv_p$-spaces. (See also Section~\ref{sec:Karami} for some comments concerning those results.) 
\end{remark}

Of course, any map $f \colon E \to E$ which is H\"older continuous with exponent $\alpha \in (0,1]$ is locally H\"older continuous in the stronger sense with the same exponent. The opposite implication, in general, does not hold.

\begin{example}
Let $\alpha \in (0,1)$ and let $E$ be a normed space. Set $f \colon E \to E$ by $f(u)=u$. Then, $\norm{f(u)-f(w)}=\norm{u-w}\leq \norm{u-w}^\alpha$ for any $u,w \in E$ with $\norm{u-w}\leq 1$. So, $f \in \Lip^\alpha_{\text{strg}}(E)$ with $L=\delta=1$. However, $f$ is not H\"older continuous on $E$ with exponent $\alpha$, since for any $u \in E\setminus \{0\}$ with $\norm{u}=1$ we have
\[
 \sup_{t > 0} \frac{\norm{f(tu)}}{\norm{tu}^\alpha} = \sup_{t>0} \, \abs{t}^{1-\alpha}=+\infty.
\]
\end{example}

Unfortunately, it is still not known whether the classes $\Lip^1(E)$ and $\Lip^1_{\text{strg}}(E)$ differ. However, if we replace the normed space $E$ with a metric one $X$, then we can already find maps $f \colon X \to X$ that are locally Lipschitz continuous in the stronger sense but are not Lipschitz continuous.

\begin{example}
Let $X:=\bigcup_{n=1}^\infty [4n,4n+1]$ be a metric subspace of $\mathbb R$ endowed with the Euclidean metric $d_{\mathbb R}$. Moreover, let $f(u)= u+4n^2$ for $u \in [4n,4n+1]$, where $n \in \mathbb N$. Note that if $u=4m+\lambda$ for some $m \in \mathbb N$ and $\lambda \in [0,1]$, then $f(u)=4m+\lambda + 4m^2=4m(m+1)+\lambda$. This means that $f(u) \in [4m(m+1),4m(m+1)+1]$. Consequently, $f$ maps $X$ into itself. It is also easy to see that $d_{\mathbb R}(f(u),f(w))=d_{\mathbb R}(u,w)$ for $u,w \in X$ with $d_{\mathbb R}(u,w)\leq 1$, as any such points must lie in the same interval $[4n,4n+1]$. However, for every $n \geq 2$ we have 
\[
 \frac{d_{\mathbb R}(f(4),f(4n))}{d_{\mathbb R}(4,4n)}=\frac{4(n-1)(n+2)}{4(n-1)}=n+2.
\]
Therefore, $f$ is not Lipschitz continuous on the whole space $X$.
\end{example}

The previous two examples ``worked'', because the considered functions were unbounded. Under the boundedness assumption, there are no differences between maps that satisfy H\"older condition and local H\"older condition in the stronger sense.

\begin{proposition}\label{prop:lip_vs_lip_strong}
Let $\alpha \in (0,1]$ and let $f \in \Lip^\alpha_{\text{strg}}(E)$. If $f$ is bounded on $E$, then $f \in \Lip^\alpha(E)$.
\end{proposition}

\begin{proof}
Since $f \in \Lip^\alpha_{\text{strg}}(E)$, there exist $L\geq 0$ and $\delta>0$ such that $\norm{f(a)-f(b)}\leq L\norm{a-b}^\alpha$ for $a,b \in E$ with $\norm{a-b}\leq \delta$. Set $M:=\sup_{v \in E}\norm{f(v)}<+\infty$ and $\Lambda:=\max\{L,2M\delta^{-\alpha}\}$. Take any two points $u,w \in E$. If $\norm{u-w}\leq \delta$, then $\norm{f(u)-f(w)}\leq L\norm{u-w}^\alpha\leq \Lambda\norm{u-w}^\alpha$. If, on the other hand, $\norm{u-w}>\delta$, then $\norm{f(u)-f(w)}\leq 2M\leq 2M\delta^{-\alpha}\norm{u-w}^{\alpha} \leq \Lambda\norm{u-w}^\alpha$. This ends the proof.
\end{proof}

A similar reasoning leads to the following result.

\begin{proposition}\label{prop:lip_strong_vs_lip_loc}
For every $\alpha \in (0,1]$ we have $\Lip^\alpha_{\text{strg}}(E) \subseteq \Lip^\alpha_{\text{bnd}}(E)$.
\end{proposition}

\begin{proof}
Let $f \in \Lip^\alpha_{\text{strg}}(E)$. Then, there exist $L\geq 0$ and $\delta>0$ such that $\norm{f(a)-f(b)}\leq L\norm{a-b}^\alpha$ for $a,b \in E$ with $\norm{a-b}\leq \delta$. We will show that $M_r:=\sup_{\norm{v}\leq r} \norm{f(v)}<+\infty$ for every $r>0$. To this end let us fix $r>0$. Moreover, set $n:=[r\delta^{-1}]$, where $[t]$ denotes the smallest integer larger than or equal to $t$. Take any $u \in E$. If $\norm{u}\leq \delta$, then $\norm{f(u)} \leq \norm{f(0)}+\norm{f(u)-f(0)}\leq \norm{f(0)}+L\delta^\alpha$. If $\delta \leq \norm{u} \leq 2\delta$, set $w:= \delta\norm{u}^{-1}u$. Then, $\norm{w}=\delta$ and $\norm{w-u}\leq \delta$. Hence, by the previous part, we obtain $\norm{f(u)}\leq \norm{f(w)}+\norm{f(u)-f(w)}\leq \norm{f(0)}+2L\delta^\alpha$. Continuing this process, we end up with
\[
 M_r \leq \sup_{\norm{v}\leq n\delta} \norm{f(v)} \leq \norm{f(0)}+nL\delta^\alpha \leq \norm{f(0)}+ (1+r\delta^{-1})L\delta^\alpha<+\infty.
\]
To end the argument it suffices now to repeat the proof of Proposition~\ref{prop:lip_vs_lip_strong} with $M$ replaced by $M_r$. 
\end{proof}

The inclusion in Proposition~\ref{prop:lip_strong_vs_lip_loc}, in general, is strict, as is shown by the following example.

\begin{example}
Let $\alpha \in (0,1]$. Consider the non-constant function $f \colon \mathbb R \to \mathbb R$ given by $f(u)=u^2$. It is clearly H\"older continuous on bounded subsets of $\mathbb R$ with exponent $\alpha$. However, for any $L,\delta>0$, if we take $u:=\frac{1}{2}L\delta^{\alpha-1}$ and $w:=\frac{1}{2}L\delta^{\alpha-1}+\delta$, we have $\abs{u-w}=\delta$. But 
\begin{align*}
\abs{f(u)-f(w)}&=\abs{u-w}\abs{u+w}=\abs{u-w}^\alpha \abs{u-w}^{1-\alpha} (L\delta^{\alpha-1}+\delta)\\
 & > \abs{u-w}^\alpha \delta^{1-\alpha} L\delta^{\alpha-1}=L\abs{u-w}^\alpha,
\end{align*}
meaning that $f \notin \Lip^\alpha_{\text{strg}}(\mathbb R)$.
\end{example}

To conclude the discussion regarding maps that exhibit local H\"older continuity in the stronger sense, let us prove their yet another equivalent characterization. 

\begin{proposition}\label{prop:eq_lip_strg}
Let $\alpha \in (0,1]$. Then, $f \in \Lip^\alpha_{\text{strg}}(E)$ if and only if for every $\eta>0$ there exists $L_\eta \geq 0$ such that $\norm{f(u)-f(w)}\leq L_\eta \norm{u-w}^\alpha$ for any $u,w \in E$ with $\norm{u-w}\leq \eta$.
\end{proposition}

\begin{proof}
Clearly, we only need to prove the necessity part. So, let us assume that there exist constants $\delta>0$ and $L\geq 0$ such that $\norm{f(a)-f(b)}\leq L\norm{a-b}^\alpha$ for $a,b \in E$ with $\norm{a-b}\leq \delta$. Also, take arbitrary $\eta>0$. Decreasing $\delta>0$, if necessary, we may assume that $k:=\eta \delta^{-1} \in \mathbb N$. Moreover, set $L_\eta:= Lk^{1-\alpha}$. Given two points $u,w \in E$ such that $\norm{u-w}\leq \eta$ let $v_j:=\frac{j}{k}w + (1-\frac{j}{k})u$ for $j=0,1,\ldots,k$. Then, $\norm{v_{j+1}-v_j}=\frac{1}{k}\norm{u-w}$. Hence, $\norm{v_{j+1}-v_j}\leq \delta$ and $\norm{f(v_{j+1})-f(v_j)}\leq L\norm{v_{j+1}-v_j}^\alpha = L k^{-\alpha}\norm{u-w}^\alpha$ for $j=0,1,\ldots,k-1$. Consequently,
\[
 \norm{f(u)-f(w)}\leq \sum_{j=0}^{k-1}\norm{f(v_{j+1})-f(v_j)} \leq Lk^{1-\alpha} \norm{u-w}^\alpha = L_\eta \norm{u-w}^\alpha.
\]
This ends the proof.
\end{proof}

\begin{remark}
Note that from Proposition~\ref{prop:eq_lip_strg} and the trivial observation that $\norm{u-w}\leq 2r$ for $u,w \in \ball_E(0,r)$ another proof of Proposition~\ref{prop:lip_strong_vs_lip_loc} follows.
\end{remark}

Clearly, for every $\alpha \in (0,1]$ we have $\Lip^\alpha_{\text{bnd}}(E) \subseteq \Lip^\alpha_{\text{comp}}(E)$ with equality if $E:=\mathbb K^n$ for $\mathbb K \in \{\mathbb R, \mathbb C\}$ and some $n \in \mathbb N$. However, if $E$ is infinite-dimensional, the notions of H\"older continuity on bounded and compact sets are not necessarily equivalent.

\begin{example}\label{ex:lip_loc_vs_lip_comp}
Set $a:=(1,\frac{1}{2^2}, \frac{1}{3^2},\ldots)$ and for $u \in c_{00}(\mathbb R)$ let $f(u)=\frac{u}{\norm{u-a}_\infty}$. Since $a\notin c_{00}(\mathbb R)$, for every $u \in c_{00}(\mathbb R)$ we have $\norm{u-a}_\infty >0$. This means that $f$ is a well-defined mapping from $c_{00}(\mathbb R)$ into itself.

Now, we are going to show that it is Lipschitz continuous on compact subsets of $c_{00}(\mathbb R)$. To this end, fix a compact set $K \subseteq c_{00}(\mathbb R)$. As the map $u \mapsto \norm{u-a}_{\infty}$ is continuous on $c_{00}(\mathbb R)$, we infer that $m:=\inf_{u \in K}\norm{u-a}_\infty >0$ and $M:=\sup_{u \in K}\norm{u-a}_\infty < +\infty$. Thus, for any $u,w \in K$, we obtain 
\begin{align*}
&\norm{f(u)-f(w)}_\infty\\
&\qquad=\norm[\Bigg]{\frac{u}{\norm{u-a}_\infty} - \frac{w}{\norm{w-a}_\infty}}_\infty\\
&\qquad=\frac{\norm[\big]{ \norm{w-a}_\infty \cdot u - \norm{u-a}_\infty \cdot w}_\infty}{\norm{u-a}_\infty \cdot \norm{w-a}_\infty}\\
&\qquad\leq \frac{\norm[\big]{ \norm{w-a}_\infty \cdot u - \norm{w-a}_\infty \cdot w}_\infty}{m^2} + \frac{\norm[\big]{ \norm{w-a}_\infty\cdot w - \norm{u-a}_\infty\cdot w}_\infty}{m^2}\\
&\qquad \leq \frac{M\norm{u-w}_\infty}{m^2} + \frac{\abs[\big]{ \norm{w-a}_\infty - \norm{u-a}_\infty}\cdot \norm{w}_\infty}{m^2}\\
&\qquad\leq \frac{M}{m^2}\norm{u-w}_\infty + \frac{M+1}{m^2}\norm{u-w}_\infty.
\end{align*}
Hence, $f$ is Lipschitz continuous on compact subsets of $c_{00}(\mathbb R)$.

Note, however, that $f$ is not locally bounded, that is, for some $r>0$ the image $f(\ball_{c_{00}}(0,r))$ is an unbounded subset of $c_{00}(\mathbb R)$. Indeed, $\norm{P_n(a)}_\infty=1$ and $\norm{P_n(a)-a}_{\infty}=\norm{Q_n(a)}_\infty=(n+1)^{-2}$ for every $n \in \mathbb N$. Hence,
\[
 \sup_{\norm{u}_\infty \leq 1} \norm{f(u)}_{\infty} \geq \sup_{n \in \mathbb N} \norm{f(P_n(a))}_{\infty}=\sup_{n \in \mathbb N} \frac{\norm{P_n(a)}_{\infty}}{\norm{Q_n(a)}_{\infty}} = \sup_{n \in \mathbb N}\ (n+1)^2=+\infty.
\] 
Therefore, $f$ cannot be Lipschitz continuous on bounded subsets of $c_{00}(\mathbb R)$. 
\end{example}

Surprisingly, even if we assume that the function $f$ is (locally) bounded, the conclusion of Example~\ref{ex:lip_loc_vs_lip_comp} remains unchanged. 

\begin{example}\label{ex:lip_loc_vs_lip_comp_2}
Let us consider the map $\psi \colon l^2(\mathbb R)\to \mathbb R$ given by
\[
 \psi(u)=\sup\dset[\big]{(2\pi+1)n \abs{u(n)}-n}{n \in \mathbb N}.
\]
Also, for any $u \in l^2(\mathbb R)$ let $f(u)=(\sin \psi(u),0,0,\ldots)$. Then, it can be shown that $f$ is a well-defined mapping acting from $l^2(\mathbb R)$ into itself (cf.~\cite{CLM}*{Example~2.160}). Furthermore, for every point $w \in l^2(\mathbb R)$ an open ball $B_{l^2}(w,r_w)$ and a constant $L_w\geq 0$ exist such that $\norm{f(u)-f(v)}_{l^2} \leq L_w\norm{u-v}_{l^2}$ for any $u,w \in B_{l^2}(w,r_w)$ -- cf.~\cite{CLM}*{pp.~63--64}. Consequently, $f$ is Lipschitz continuous on compact subsets of $l^2(\mathbb R)$ (see~\cite{cobzas_book}*{Theorem~2.1.6 and Remark~2.1.7}).  It is also bounded on $l^2(\mathbb R)$, that is, $\sup_{u\in l^2(\mathbb R)} \norm{f(u)}_{l^2}<+\infty$.

Now, we are going to show that $f$ is not Lipschitz continuous on bounded subsets of $l^2(\mathbb R)$. By $e_k$ let us denote the $k$-th unit vector of $l^2(\mathbb R)$, that is, the vector which consists of all zeros and a one on $k$-th position. Note that 
$\psi(e_k)=2\pi k$ and $\psi\bigl((1-k^{-2})e_k\bigr)=2\pi k - (2\pi+1)k^{-1}$. And thus, $f(e_k)=0$ and $f\bigl((1-k^{-2})e_k\bigr)=(-\sin[(2\pi+1)k^{-1}],0,0,\ldots)$. This, in turn, implies that for $k \geq 5$ we have
\[
 \frac{\norm[\big]{f(e_k) - f\bigl((1-\tfrac{1}{k^{2}})e_k\bigr)}_{l^2}}{\norm[\big]{e_k - (1-\tfrac{1}{k^{2}})e_k}_{l^2}} =  k^2 \norm[\big]{f(e_k) - f\bigl((1-\tfrac{1}{k^{2}})e_k\bigr)}_{l^2} =   k^2 \cdot \sin \tfrac{2\pi+1}{k} \geq k^2 \cdot \tfrac{2}{\pi} \cdot \tfrac{2\pi+1}{k} \geq 4k.
\]
So, $f$ is not Lipschitz continuous on bounded subsets of $l^2(\mathbb R)$.
\end{example}

\section{Acting conditions}
\label{sec:action_conditions}

In the theory of composition operators the problem that is the most fundamental and comes before all the others is the problem of characterization of \emph{acting conditions}. In other words, given two sequence spaces $X$ and $Y$, we want to find all the mappings $f$ that generate composition operators $C_f$ acting from $X$ into $Y$.

For simplicity's sake, hereafter let us adopt the following convention. Whenever we will be talking about the composition operator $C_f$, we will always assume that it is generated by the map $f \colon E \to E$, where $E$ is a normed space (over the field of either real or complex numbers). 

\subsection{Composition operators acting into $bv_p(E)$}
\label{sec:AC_into_bv}

As the title suggests, in this subsection we will be interested in characterizing acting conditions for those composition operators that act into $bv_p(E)$. Let us begin with a lemma, which may be of independent interest.

\begin{lemma}\label{lem:into_bv_q_f_continous}
Let $p,q \in [1,+\infty)$. Moreover, let $X \in \{bv_p(E), c(E), l^\infty(E)\}$ and $Y \in \{l^p(E),c_0(E)\}$. 
\begin{enumerate}[label=\textup{(\alph*)}]
 \item\label{it:l1a} If $C_f$ maps $X$ into $bv_q(E)$, then $f$ is continuous on $E$.
 
 \item\label{it:l1b} If $C_f$ maps $Y$ into $bv_q(E)$, then $f$ is continuous at $0$.
\end{enumerate}
\end{lemma}

\begin{proof}
We begin with the proof of part~\ref{it:l1a}. Suppose on the contrary that the composition operator $C_f$ maps $X$ into $bv_q(E)$, but $f$ is not continuous at a point $u \in E$.  Then, there exist a number $\varepsilon>0$ and a sequence $(u_n)_{n \in \mathbb N}$ of elements of $E$ such that $\lim_{n \to \infty} u_n = u $ and $\norm{f(u_n)-f(u)}\geq \varepsilon$ for every $n \in \mathbb N$. Of course, without loss of generality, we may assume that $\norm{u_n - u}\leq \frac{1}{n^2}$ for $n \in \mathbb N$. Now, we define $x \colon \mathbb N \to E$ by the formula
\[
 x(n):=\begin{cases}
        u_{\frac{1}{2}n} & \text{for $n \in \mathbb N$ even,}\\
				u & \text{for $n \in \mathbb N$ odd.}
				\end{cases}				
\]
It is easy to check that $x \in X$. On the other hand, 
\[
\sum_{n=1}^\infty \norm{f(x(n+1))-f(x(n))}^q \geq \sum_{n=1}^\infty \norm{f(u_n)-f(u)}^q \geq \sum_{n=1}^\infty \varepsilon^q = +\infty.
\]
Hence, $C_f(x)\notin bv_q(E)$. This shows that $f$ must be continuous on $E$.

The proof of part~\ref{it:l1b} is similar to the proof of~\ref{it:l1a}. Moreover, it suffices to consider the case $Y:=l^1(E)$ only. Therefore, we will skip it.
\end{proof}

Part~\ref{it:l1b} of Lemma~\ref{lem:into_bv_q_f_continous} in the case $X:=l^p(E)$ cannot be strengthened. As the following example shows there are maps that are continuous only at zero and generate composition operators acting from $l^p(E)$ into $bv_q(E)$. On the other hand, in Theorem~\ref{thm:c_0_into_bv} below we will see that for $X:=c_0(E)$ a result much stronger than  Lemma~\ref{lem:into_bv_q_f_continous}~\ref{it:l1b} is true.

\begin{example}\label{ex:lp_into_bvq}
Let $p,q \in [1,+\infty)$ and let $f \colon \mathbb R \to \mathbb R$ be given by $f(u)=\abs{u}^{\frac{p}{q}}\cdot \chi_{\mathbb Q}(u)$. It is clear that the function $f$ is continuous only at zero. Furthermore, for any $x \in l^p(E)$ and $n \in \mathbb N$ we have $\abs{f(x(n+1))-f(x(n))}^q \leq 2^q \abs{x(n+1)}^p + 2^q \abs{x(n)}^p$. Thus, $C_f(l^p(\mathbb R)) \subseteq bv_q(\mathbb R)$. 
\end{example}

Before we will be ready to discuss acting conditions for $C_f$, we need yet another technical lemma. The idea behind this result and its proof comes from Dedagich and Zabrejko (see~\cite{DZ}*{p.~88}). 

\begin{lemma}\label{lem:1_bv_version}
Let $p,q \in [1,+\infty)$ and let $X \in \{c_0(E), l^p(E)\}$. Moreover, assume that $f(0)=0$. If $C_f$ maps $X$ into $bv_q(E)$, then for every $\varepsilon>0$ there exist two numbers $\delta>0$ and $N \in \mathbb N$ such that $(C_f\circ Q_N)(\ball_X(0,\delta))\subseteq \ball_{bv_q}(0,\varepsilon)$.
\end{lemma}

\begin{proof}
Suppose on the contrary that there is $\varepsilon>0$ such that for each $n \in \mathbb N$ we can find an element $x_n \in \ball_X(0,2^{-n})$ with $\norm{(C_f\circ Q_n)(x_n)}_q > \varepsilon$. Then, keeping $n$ fixed, by Lemma~\ref{lem:into_bv_q_f_continous}~\ref{it:l1b} we have $\lim_{m \to \infty} f(x_n(m))=0$. This, in turn, implies that 
\[
 \lim_{m \to \infty} \norm{C_f\circ (Q_n - Q_m)(x_n)}_{q} = \norm{(C_f \circ Q_n)(x_n)}_{q}>\varepsilon.
\]
Hence, for each $n \in \mathbb N$ there is $m_n > n$ such that $\norm{C_f\circ (Q_n - Q_{m_n})(x_n)}_{q}>\varepsilon$. Thus, starting with $n_1:=1$ we can construct inductively a sequence of indices setting $n_{k+1}:=m_{n_k}$, where $k \in \mathbb N$, for which the following estimates hold: $\norm{(Q_{n_k}-Q_{n_{k+1}})(x_{n_k})}_X \leq \norm{x_{n_k}}_X \leq 2^{-k}$ and $\norm{C_f\circ (Q_{n_k} - Q_{n_{k+1}})(x_{n_k})}_{q}>\varepsilon$. Here, by $\norm{\cdot}_X$ we denote the default norm on the space $X$. Now, let $x:=(0,x_{n_1}(n_1+1),\ldots,x_{n_1}(n_2),0,x_{n_2}(n_2+1),\ldots,x_{n_2}(n_3),0,\ldots)$. It is easy to check that $x \in X$. However,
\begin{align*}
&\sum_{j=1}^{\infty} \norm[\big]{C_f(x)(j+1)-C_f(x)(j)}^q\\
&\quad = \sum_{k=1}^\infty \Biggl(\norm[\big]{f(x_{n_k}(n_k+1))}^q + \sum_{j=n_{k}+1}^{n_{k+1}-1} \norm[\big]{f(x_{n_k}(j+1))-f(x_{n_k}(j))}^q
 + \norm[\big]{f(x_{n_k}(n_{k+1}))}^q\Biggr)\\
&\quad = \sum_{k=1}^\infty \norm[\big]{C_f\circ (Q_{n_k} - Q_{n_{k+1}})(x_{n_k})}_{q}^q \geq  \sum_{k=1}^\infty\varepsilon^q =+\infty.
\end{align*}
This means that $C_f(x)\notin bv_q(E)$. The proof is complete.
\end{proof}

Finally, we are in position to study acting conditions for the composition operator $C_f$ with $bv_p(E)$ as the target space. We begin with operators acting from $c_0(E)$. Then, we move to the spaces $c(E)$ and $l^p(E)$ for $p \in [1,+\infty]$. At the end, we will study composition operators acting between $bv_p$-spaces.  

\begin{theorem}\label{thm:c_0_into_bv}
Let $p \in [1,+\infty)$. The composition operator $C_f$ maps $c_0(E)$ into $bv_p(E)$ if and only if its generator $f$ is locally constant at $0$, that is, there exists a number $\delta>0$ such that $f|_{\ball_E(0,\delta)}$ is constant.
\end{theorem}

\begin{proof}
First, assume that the map $f$ is constant on some closed ball $\ball_E(0,\delta)$. Take $x \in c_0(E)$ and choose an index $N \in \mathbb N$ so that $\norm{x(n)} \leq \delta$ for $n \geq N$. Then, for all $n \geq N$ we have $C_f(x)(n+1)=C_f(x)(n)$. This implies that $C_f(x)\in bv_p(E)$.

Now, let us assume that $C_f$ maps $c_0(E)$ into $bv_p(E)$. Define the map $g \colon E \to E$ by setting $g(v):=f(v)-f(0)$. Then, clearly, $g(0)=0$. Furthermore, the composition operator $C_g$ maps $c_0(E)$ into $bv_p(E)$. Hence, in view of Lemma~\ref{lem:1_bv_version}, there exist two numbers $\delta >0$ and $N \in \mathbb N$ such that $\norm{(C_g\circ Q_N)(x)}_{p}\leq 1$ for all $\norm{x}_\infty \leq \delta$. 

Suppose that there are two points $u,w \in \ball_E(0,\delta)$ such that $f(u)\neq f(w)$. Also, take $m \in \mathbb N$ so that $m^{\frac{1}{p}}\cdot \norm{f(u)-f(w)}>1$ and define $y \in \ball_{c_0}(0,\delta)$ by the formula 
\[
 y(n):=\begin{cases}
        0, & \text{if $1\leq n \leq N$, or $n > N+2m$,}\\
				u, & \text{if $N+1 \leq n \leq N+2m$ is odd,}\\ 
				w, & \text{if $N+1 \leq n \leq N+2m$ is even.}
				\end{cases}
\]
Then, $\norm{(C_g\circ Q_N)(y)}_p\geq m^{\frac{1}{p}} \cdot \norm{f(u)-f(w)}>1$, which is impossible. Thus, $f$ is locally constant at $0$.
\end{proof}

Surprisingly, from the above theorem with almost no effort we obtain a result for $c(E)$ and $l^{\infty}(E)$. 

\begin{theorem}\label{thm:c_l_infinity_into_bv}
Let $p \in [1,+\infty)$ and let $X\in \{c(E),l^\infty(E)\}$. The composition operator $C_f$ maps $X$ into $bv_p(E)$ if and only if $f$ is a constant function.
\end{theorem}

\begin{proof}
Clearly, we need to show only the necessity part. 

First, note that since the space $X$ contains all the constant sequences, the function $f$ generates a composition operator that maps $X$ into $bv_p(E)$ if and only if its translate $t \mapsto f(t+e)$ by a vector $e \in E$ does. Now, let us take any two distinct points $u,w \in E$. Also, let us set $K:=\dset{tu+(1-t)w}{t \in [0,1]}$. By Theorem~\ref{thm:c_0_into_bv}, the fact that $c_0(E)\subseteq X$ and our initial observation, for each $v \in K$ there exists $\delta_v>0$ such that $f$ is constant on $\ball_E(v,\delta_v)$. In view of the compactness of $K$, this implies that we can write $K=\bigcup_{j=1}^n (K\cap B_E(v_j,\delta_{v_j}))$ for some finite family $v_1,\ldots,v_n$ of points in $K$ and the corresponding family $\delta_{v_1},\ldots,\delta_{v_n}$ of radii. Note that each of the sets $K\cap B_E(v_j,\delta_{v_j})$ has a non-empty intersection with at least one of the remaining ones. Otherwise, $K$ would be disconnected, which is absurd. Thus, inductively, we can prove that $f$ is constant on the whole set $K$. In particular, $f(u)=f(w)$. As the points $u,w$ were arbitrary, this shows that $f \colon E \to E$ is a constant map.
\end{proof}

As the next step of our investigation let us discuss composition operators having $l^p$-spaces as their domains. 

\begin{theorem}\label{thm:lp_into_bvq}
For any $p,q \in [1,+\infty)$ the following conditions are equivalent\textup:
\begin{enumerate}[label=\textup{(\roman*)}] 
 \item\label{it:lp_into_bv_i} the composition operator $C_f$ maps $l^p(E)$ into $bv_q(E)$, 
 \item\label{it:lp_into_bv_ii} there exist $\delta>0$ and $L\geq 0$ such that $\norm{f(u)-f(w)}\leq L\norm{u}^{\frac{p}{q}}+L\norm{w}^{\frac{p}{q}}$ for all $u,w \in \ball_E(0,\delta)$,
 \item\label{it:lp_into_bv_iii} there exist $\delta>0$ and $L\geq 0$ such that $\norm{f(u)-f(0)}\leq L\norm{u}^{\frac{p}{q}}$ for every $u \in \ball_E(0,\delta)$.
\end{enumerate}
\end{theorem}

The proof of the equivalence of the first two conditions of Theorem~\ref{thm:lp_into_bvq} is inspired by the celebrated result due to Josephy, which characterizes acting conditions for a composition operator in the space $BV_1([a,b], \mathbb R)$ (see, e.g.,~\cite{ABM}*{Theorem~5.9},~\cite{BBK}*{Theorem~6.4.13} or \cite{J}*{Theorem~4}). 

\begin{proof}
The proof of the implication $\ref{it:lp_into_bv_ii} \Rightarrow \ref{it:lp_into_bv_i}$ is obvious (cf.~Example~\ref{ex:lp_into_bvq}). So without further ado, let us proceed to the opposite one $\ref{it:lp_into_bv_i} \Rightarrow \ref{it:lp_into_bv_ii}$. We will argue by contradiction. By Lemma~\ref{lem:into_bv_q_f_continous}~\ref{it:l1b} there exists a radius $\delta>0$ such that $\norm{f(0)-f(u)}\leq 1$ for all $u \in \ball_E(0,\delta)$. Set $M:=1+\sup_{\norm{u}\leq \delta}\norm{f(u)}$. Also, suppose that for any $n \in \mathbb N$ we can find two points $u_n, w_n \in \ball_E(0,\delta)$ with
\[
 \norm{f(u_n)-f(w_n)} > 2^{1+\frac{1}{q}}Mn^2\bigl(\norm{u_n}^{\frac{p}{q}} + \norm{w_n}^{\frac{p}{q}}\bigr).
\] 
Note that $u_n$ and $w_n$ cannot be zero at the same time. Set
\[ 
 m_n:=\bigl[n^{-2q}\bigl(\norm{u_n}^{p} +\norm{w_n}^{p}\bigr)^{-1}\bigr],
\]
where $[t]$ denotes the greatest integer smaller than or equal to $t$. Observe that $m_n$ is well-defined and, moreover, $m_n \geq 2$ for $n \in \mathbb N$. Indeed, we have $2M \geq  2^{1+\frac{1}{q}}Mn^2\bigl(\norm{u_n}^{\frac{p}{q}} + \norm{w_n}^{\frac{p}{q}}\bigr)$. This, in turn implies, that
\[
 n^{-2q} \geq 2 \bigl(\norm{u_n}^{\frac{p}{q}} + \norm{w_n}^{\frac{p}{q}}\bigr)^q \geq 2\bigl(\norm{u_n}^p + \norm{w_n}^p\bigr),
\]
whence the claim follows. Therefore, 
\[
 \tfrac{1}{2}n^{-2q}\bigl(\norm{u_n}^{p} +\norm{w_n}^{p}\bigr)^{-1}\leq m_n \leq n^{-2q}\bigl(\norm{u_n}^{p} +\norm{w_n}^{p}\bigr)^{-1}
\]
for $n \in \mathbb N$. Consider now the following sequence $x$:
\[
 \underbrace{u_1, w_1, u_1, w_1, \ldots, u_1,w_1,}_{\text{$2m_1$ elements}}\ \underbrace{u_2, w_2, u_2, w_2, \ldots, u_2,w_2,}_{\text{$2m_2$ elements}}\ \ldots
\]
Then,  
\begin{align*}
 \sum_{n=1}^\infty \norm{x(n)}^p& = \sum_{n=1}^\infty m_n\norm{u_n}^p +\sum_{n=1}^\infty m_n\norm{w_n}^p \leq \sum_{n=1}^\infty \frac{1}{n^{2q}} < +\infty.
\end{align*}
Thus, $x \in l^p(E)$. However, on the other hand, we have
\begin{align*}
 \sum_{n=1}^\infty \norm[\big]{f(x(n+1))-f(x(n))}^q &\geq \sum_{n=1}^\infty m_n \norm{f(u_n)-f(w_n)}^q\\
& \geq \sum_{n=1}^\infty 2^{q+1}M^q m_n n^{2q}\bigl(\norm{u_n}^{\frac{p}{q}}+\norm{w_n}^{\frac{p}{q}}\bigr)^q\\
& \geq \sum_{n=1}^\infty 2^{q+1}M^q m_n n^{2q}\bigl(\norm{u_n}^p+\norm{w_n}^p\bigr)\\
& \geq \sum_{n=1}^\infty (2M)^q = +\infty. 
\end{align*}
Hence, $C_f(x)\notin bv_q(E)$ -- a contradiction.

The equivalence of~\ref{it:lp_into_bv_ii} and~\ref{it:lp_into_bv_iii} follows easily from the triangle inequality.
\end{proof}

Lastly, we turn our attention to composition operators acting between $bv_p$-spaces.  Unsurprisingly, this will be the most challenging part of this subsection due to the distinct nature of $bv_1(E)$ and $bv_p(E)$ for $p>1$. We start with a general result, the proof of which once again is based on Josephy's approach.

\begin{proposition}\label{prop:josephy_bv}
Let $p,q \in [1,+\infty)$. If the composition operator $C_f$ maps $bv_p(E)$ into $bv_q(E)$, then $f$ is H\"older continuous on compact subsets of $E$ with exponent $\frac{p}{q}$. 
\end{proposition}

\begin{proof}
Let us assume that $C_f$ maps $bv_p(E)$ into $bv_q(E)$. Also, suppose that there is a compact subsets $K$ of $E$ so that  for any $n \in \mathbb N$ we can find two points $u_n, w_n \in K$ with
\[
 \norm{f(u_n)-f(w_n)} > 2^{1+\frac{1}{q}}Mn^2\norm{u_n-w_n}^{\frac{p}{q}},
\]
where $M:=1+\sup_{v \in K}\norm{f(v)}$; note that by Lemma~\ref{lem:into_bv_q_f_continous}~\ref{it:l1a} the map $f$ is continuous on $E$, and hence $M<+\infty$.
The points $u_n$, $w_n$ must be distinct, and moreover $\norm{u_n-w_n}^p\leq \frac{1}{2}n^{-2q}$ for $n \in \mathbb N$.  Passing to a subsequence if necessary, we may assume that the sequence $(u_n)_{n \in \mathbb N}$ is convergent to a point $u_\ast \in K$ and $\norm{u_n-u_\ast}^p\leq \frac{1}{2}n^{-2q}$ for $n \in \mathbb N$. Then, also $w_n \to u_\ast$ as $n \to +\infty$, and $\norm{w_n -u_{n+1}}^p\leq 3^{p+1}n^{-2q}$ for $n \in \mathbb N$. Now, let $m_n:=[n^{-2q}\norm{u_n-w_n}^{-p}]$. As in the proof of Theorem~\ref{thm:lp_into_bvq}, $[t]$ denotes the greatest integer smaller than or equal to $t$. Observe that $m_n \geq 2$ for every $n \in \mathbb N$. Hence, $\frac{1}{2}n^{-2q}\norm{u_n-w_n}^{-p}\leq m_n \leq n^{-2q}\norm{u_n-w_n}^{-p}$ for $n \in \mathbb N$. Consider now the following sequence $x$:
\[
 \underbrace{u_1, w_1, u_1, w_1, \ldots, u_1,w_1,}_{\text{$2m_1$ elements}}\ \underbrace{u_2, w_2, u_2, w_2, \ldots, u_2,w_2,}_{\text{$2m_2$ elements}}\ \ldots
\]
Then,     
\begin{align*}
 \sum_{n=1}^\infty \norm{x(n+1)-x(n)}^p&= \sum_{n=1}^\infty (2m_n-1)\norm{u_n-w_n}^p + \sum_{n=1}^\infty \norm{w_n - u_{n+1}}^p\\
 & \leq \sum_{n=1}^\infty \frac{2}{n^{2q}} + \sum_{n=1}^\infty \frac{3^{p+1}}{n^{2q}} < +\infty.
\end{align*}
This shows that $x$ belongs to $bv_p(E)$. Similarly to what we did in the proof of Theorem~\ref{thm:lp_into_bvq} we can also show that
\begin{align*}
\sum_{n=1}^\infty \norm[\big]{f(x(n+1))-f(x(n))}^q &\geq \sum_{n=1}^\infty (2M)^q = +\infty. 
\end{align*} 
Hence, $C_f(x) \notin bv_q(E)$ -- a contradiction. Thus, $f$ is H\"older continuous on compact subsets of $E$ with exponent~$\frac{p}{q}$.
\end{proof}

From Remark~\ref{rem:holder_alpha_greater_than_1} and Proposition~\ref{prop:josephy_bv} we immediately obtain the following corollary.

\begin{corollary}\label{cor:bv_p_bv_q}
Let $1\leq q<p < +\infty$. The composition operator $C_f$ acts between $bv_p(E)$ and $bv_q(E)$ if and only if $f$ is a constant map. 
\end{corollary}

Another corollary to Proposition~\ref{prop:josephy_bv} deals with composition operators mapping $bv_1(E)$ into $bv_p(E)$. 

\begin{theorem}\label{thm:1}
Let $p \in [1,+\infty)$ and let $E$ be a Banach space. The composition operator $C_f$ maps $bv_1(E)$ into $bv_p(E)$ if and only if $f$ is H\"older continuous on compact subsets of $E$ with exponent~$\frac{1}{p}$.
\end{theorem}

\begin{proof}
Thanks to Proposition~\ref{prop:josephy_bv} we need to show the sufficiency part only. To this end let us fix a sequence $x \in bv_1(E)$. As $E$ is a Banach space, the sequence $x$ is convergent. Denote its limit in $E$ by $x_\ast$. Moreover, let $K:=\{x_\ast\}\cup \dset{x(n) \in E}{n \in \mathbb N}$. Clearly, $K$ is a compact subset of $E$. Hence,
\[
 \sum_{n=1}^\infty \norm[\big]{f(x(n+1))-f(x(n))}^p \leq L_K^p \sum_{n=1}^\infty \norm{x(n+1)-x(n)}<+\infty;
\]
here, $L_K$ is the H\"older constant of $f|_K$. This shows that $C_f(x)\in bv_p(E)$.
\end{proof}

A very natural question concerning Theorem~\ref{thm:1} is whether it also holds if $E$ is not complete. As the following example shows, the answer is negative.   

\begin{example}\label{ex:completeness_essential}
Let $a:=(1,\frac{1}{2^2}, \frac{1}{3^2},\ldots)$ and for $u \in c_{00}(\mathbb R)$ let $f(u)=\frac{u}{\norm{u-a}_\infty}$. From Example~\ref{ex:lip_loc_vs_lip_comp} we know that $f$ is Lipschitz continuous on compact subsets of $c_{00}(\mathbb R)$. However, as we are going to show, it does not generate a composition operator acting between $bv_1(c_{00}(\mathbb R))$ and $bv_1(c_{00}(\mathbb R))$. Define $x \colon  \mathbb N \to c_{00}(\mathbb R)$ by $x(n):=P_n(a)$. Since $\norm{x(n+1)-x(n)}_\infty=(n+1)^{-2}$ for $n \in \mathbb N$, we see that $x \in bv_1(c_{00}(\mathbb R))$. But,
\begin{align*}
 \sum_{n=1}^\infty \norm[\big]{f(x(n+1))-f(x(n))}_\infty &= \sum_{n=1}^\infty \norm[\big]{f(P_{n+1}(a))-f(P_n(a))}_\infty\\
&=\sum_{n=1}^\infty \sup_{k \in \mathbb N} \abs[\Bigg]{ \frac{P_{n+1}(a)(k)}{\norm{Q_{n+1}(a)}_\infty}- \frac{P_{n}(a)(k)}{\norm{Q_n(a)}_\infty}}\\
&\geq \sum_{n=1}^\infty \abs[\Bigg]{ \frac{P_{n+1}(a)(n+1)}{\norm{Q_{n+1}(a)}_\infty}- \frac{P_{n}(a)(n+1)}{\norm{Q_n(a)}_\infty}}\\
&= \sum_{n=1}^\infty  \frac{(n+2)^2}{(n+1)^2}=+\infty.
\end{align*}
Hence, $C_f(x) \notin bv_1(c_{00}(\mathbb R))$. This highlights the significance of the completeness of the normed space $E$ in Theorem~\ref{thm:1}.
\end{example}

So far, we have characterized acting conditions for $C_f \colon bv_1(E) \to bv_p(E)$ for any $p\in [1,+\infty)$ and $C_f \colon bv_p(E) \to bv_q(E)$ for $1\leq q < p <+\infty$. Now, we will investigate composition operators mapping $bv_p(E)$ into $bv_q(E)$ for $1<p\leq q < +\infty$. In a parallel theory of spaces of functions of bounded Wiener variation, it is well-known that a composition operator maps $BV_p([a,b],\mathbb R)$ into $BV_q([a,b],\mathbb R)$ for $1 \leq p \leq q <+\infty$ if and only if its generator is H\"older continuous on bounded subsets of $\mathbb R$ with exponent $\frac{p}{q}$ (see, e.g.,~\cite{ABM}*{Theorem~5.12}). Thus, it would be natural to expect a similar result in the case of sequence $bv_p$-spaces, possibly with bounded subsets of $E$ replaced by compact ones. Surprisingly, as the following example shows, this is not true.

\begin{example}
Let $f \colon \mathbb R \to \mathbb R$ be given by the formula $f(u)=u^2$. Clearly, the function $f$ is Lipschitz continuous on bounded/compact subsets of $\mathbb R$. Now, let $x \colon \mathbb N \to \mathbb R$ be defined by $x(n):=\sum_{k=1}^n k^{-\frac{3}{4}}$. It is evident that $x$ is an element of $bv_2(\mathbb R)$. However, the sequence $n \mapsto f(x(n))$ is not. Indeed, for any $n \geq 15$ we have
\begin{align*}
 \sum_{k=1}^n k^{-\frac{3}{4}} & \geq \int_1^{n+1} x^{-\frac{3}{4}} dx = 4(n+1)^{\frac{1}{4}}-4 = (n+1)^{\frac{1}{4}} + 3(n+1)^{\frac{1}{4}}-4\\
& \geq (n+1)^{\frac{1}{4}} + 3\cdot 16^{\frac{1}{4}} - 4 = (n+1)^{\frac{1}{4}} + 2>(n+1)^{\frac{1}{4}}\geq n^{\frac{1}{4}}.
\end{align*}
And hence,
\begin{align*}
 \abs[\big]{f(x(n+1))-f(x(n))}^2 &= \abs[\big]{x(n+1)-x(n)}^2 \cdot \abs[\big]{x(n+1)+x(n)}^2\\
& = (n+1)^{-\frac{3}{2}} \cdot \Biggl(\sum_{k=1}^n 2k^{-\frac{3}{4}} + (n+1)^{-\frac{3}{4}}\Biggr)^2\\
& \geq 4(n+1)^{-\frac{3}{2}} \cdot \Biggl(\sum_{k=1}^n k^{-\frac{3}{4}}\Biggr)^2 \geq 4(2n)^{-\frac{3}{2}} \cdot n^{\frac{1}{2}} = \sqrt{2}n^{-1}.
\end{align*}
Thus, $f$ does not generate a composition operator that acts between $bv_2(\mathbb R)$ and $bv_2(\mathbb R)$.
\end{example}

To obtain the characterization of acting conditions for $C_f \colon bv_p(E) \to bv_q(E)$, where $1<p\leq q<+\infty$, we need to use the local H\"older condition in the stronger sense.

\begin{theorem}\label{thm:bv_p_int_bv_q_strong}
Let $1<p\leq q < +\infty$. The composition operator $C_f$ maps $bv_p(E)$ into $bv_q(E)$ if and only if $f$ is locally H\"older continuous in the stronger sense with exponent $\frac{p}{q}$.
\end{theorem}

\begin{proof}
We begin with the proof of the sufficiency part. Assume that $f$ is locally H\"older continuous in the stronger sense, that is, there exist $\delta>0$ and $L\geq 0$ such that $\norm{f(u)-f(w)}\leq L\norm{u-w}^{\frac{p}{q}}$ for $u,w \in E$ with $\norm{u-w}\leq \delta$. Fix $x \in bv_p(E)$. Then, there exists an index $N \in \mathbb N$ such that $\norm{x(n+1)-x(n)}\leq \delta$ for all $n \geq N$. This, in turn, implies that
\[
 \sum_{n=N}^\infty \norm{f(x(n+1))-f(x(n))}^q \leq L^q \sum_{n=N}^\infty \norm{x(n+1)-x(n)}^p < +\infty.
\]
Hence, $C_f(x) \in bv_q(E)$.

Now, let us proceed to the second part of the demonstration. It will be similar to the proof of Proposition~\ref{prop:josephy_bv}. This time, however, because we have no compact set at our disposal and  for $p>1$ the space $bv_p(E)$ is not included in $c(E)$, the definition of the sequence $x$ will be slightly more complex. So, let us assume that $C_f$ maps $bv_p(E)$ into $bv_q(E)$. Also, suppose that for every $n \in \mathbb N$ there exist two distinct points $u_n,w_n \in E$ such that
\[
 \norm{u_n-w_n}\leq 2^{-\frac{1}{p}}n^{-\frac{2q}{p}} \quad \text{and} \quad \norm{f(u_n)-f(w_n)}> 2^{\frac{1}{q}}n^2\norm{u_n-w_n}^{\frac{p}{q}}.
\]
Then, $n^{-2q}\norm{u_n-w_n}^{-p} \geq 2$. Hence, setting $m_n:=[n^{-2q}\norm{u_n-w_n}^{-p}]$, we have $\frac{1}{2}n^{-2q}\norm{u_n-w_n}^{-p} \leq m_n \leq n^{-2q}\norm{u_n-w_n}^{-p}$ for $n \in \mathbb N$; here, as before, $[t]$ stands for the greatest integer smaller than or equal to $t$. Further, keeping $n$ fixed, let $k_n \in \mathbb N$ be so that $\norm{u_{n+1}-w_n}^p \cdot k_n^{1-p} \leq 2^{-n}$. And, for any $j=0,\ldots,k_n$ let $v^n_j:=\frac{j}{k_n}u_{n+1} + (1-\frac{j}{k_n})w_n$. Now, define $x \colon \mathbb N \to E$ as
\[
  \underbrace{u_1, w_1, u_1, w_1, \ldots, u_1,w_1,}_{\text{$2m_1$ elements}}\ v_1^1,v_2^1,\ldots,v_{k_1-1}^1,\ \underbrace{u_2, w_2, u_2, w_2, \ldots, u_2,w_2,}_{\text{$2m_2$ elements}}\ v_1^2,v_2^2,\ldots,v_{k_2-1}^2,\ldots 
\]
Observe that
\begin{align*}
 \sum_{n=1}^\infty \norm{x(n+1)-x(n)}^p & = \sum_{n=1}^\infty (2m_n-1)\norm{u_n-w_n}^p + \sum_{n=1}^\infty \sum_{j=0}^{k_n-1} \norm{v_{j+1}^n-v_j^n}^p\\
& \leq \sum_{n=1}^\infty 2m_n \norm{u_n-w_n}^p + \sum_{n=1}^\infty \norm{u_{n+1}-w_n}^p \cdot k_n^{1-p}\\
& \leq \sum_{n=1}^\infty \frac{2}{n^{2q}} + \sum_{n=1}^\infty \frac{1}{2^n}<+\infty.
\end{align*}
Thus, $x \in bv_p(E)$. Using a similar reasoning to the one used in the proof of Proposition~\ref{prop:josephy_bv}, we can show that $C_f(x)\notin bv_q(E)$. Therefore, $f$ is locally H\"older continuous in the stronger sense with exponent $\frac{p}{q}$. 
\end{proof}

\begin{remark}\label{rem:producing_sequences_in_bv_p}
One tool we used in the above proof is especially interesting and worth highlighting. Given a sequence $(u_n)_{n \in \mathbb N}$ of elements of a normed space $E$ and a parameter $p>1$,  we can always construct $x \in bv_p(E)$ which contains $(u_n)_{n \in \mathbb N}$ as its subsequence.

The idea behind this construction can be summarized as follows. We connect the consecutive terms $u_n$ and $u_{n+1}$ with the segment $[u_n,u_{n+1}]:=\dset{tu_{n+1}+(1-t)u_n}{t \in [0,1]}$. Then, we divide $[u_n,u_{n+1}]$ into $k_n$ subsegments $[v_{j}^n,v_{j+1}^n]$, where $j=0,\ldots,k_n-1$ and $v_j^n:=\frac{j}{k_n}u_{n+1} + (1-\frac{j}{k_n})u_n$. The numbers $k_n$ are chosen so that the series $\sum_{n=1}^\infty \norm{u_{n+1}-u_n}^p \cdot k_n^{1-p}$ is convergent. For example, we can assume that $k_n$ is such that $\norm{u_{n+1}-u_n}^p \cdot k_n^{1-p} \leq 2^{-n}$ for $n \in \mathbb N$. (This is exactly what we did in the proof of Theorem~\ref{thm:bv_p_int_bv_q_strong}.) Then, as the sequence $x$, we take all the endpoints of the subsegments $[v_{j}^n,v_{j+1}^n]$ arranged in an ``ascending'' order, that is, $x:=(v_0^1,v_1^1,\ldots,v_{k_1}^1,v_1^2,v_2^2,\ldots,v_{k_2}^2,v_1^3,v_2^3,\ldots)$.

Note that this construction does not work for $p=1$, even when $E:=\mathbb R$. Indeed, the sequence of all positive integers $(1,2,3,\ldots)$ cannot be a subsequence of any $x \in bv_1(\mathbb R)$, because every sequence of bounded variation (with $p=1$) is bounded.
\end{remark}

\subsection{Composition operators acting from $bv_p(E)$}

In this subsection we will focus on characterizing acting conditions for those composition operators that act from the space $bv_p(E)$. We will start with the following very simple result. We call it a ``theorem'' only because it provides necessary and sufficient conditions for $C_f$ to map $bv_p(E)$ into $c_0(E)$ and $l^q(E)$.

\begin{theorem}\label{thm:bv_p_into_c_l}
Let $p,q \in [1,+\infty)$ and let $X \in \{c_0(E),l^q(E)\}$. The composition operator $C_f$ maps $bv_p(E)$ into $X$ if and only if its generator $f$ is the zero function.
\end{theorem}

\begin{proof}
Clearly, it suffices to consider the case $X=c_0(E)$ only. Suppose that $f(u)\neq 0$ for some $u \in E$. Then, the constant sequence $x:=(u,u,u,\ldots)$ belongs to $bv_p(E)$ for any $p\geq 1$. However, $C_f(x)$ cannot belong to $c_0(E)$, as this space does not contain non-zero constant sequences.
\end{proof}

Now, we will move to the case when the target space is $l^\infty(E)$.

\begin{theorem}\label{thm:bv_1_into_l}
Let $E$ be a Banach space. The composition operator $C_f$ maps $bv_1(E)$ into $l^\infty(E)$ if and only if $f$ is bounded on compact subsets of $E$.
\end{theorem}

\begin{proof}
First, let us assume that $f$ is bounded on compact subsets of $E$. Take $x \in bv_1(E)$. Then, $x \in c(E)$. Consequently, $K:=\{x_\ast\} \cup \dset{x(n) \in E}{n \in \mathbb N}$ is a compact subset of $E$; here $x_\ast$ denotes the limit of the sequence $x$. And so, there exists a constant $r>0$ such that $\norm{f(x(n))} \leq r$ for $n \in \mathbb N$. Thus, $C_f(x) \in l^\infty(E)$.

Now, let as assume that $C_f$ maps $bv_1(E)$ into $l^\infty(E)$, and suppose that for some compact subset $K$ of $E$ we have $\sup_{u \in K}\norm{f(u)}=+\infty$. Then, there is a sequence $(u_n)_{n \in \mathbb N}$ of elements of $K$ such that $\norm{f(u_n)} \to +\infty$ as $n \to +\infty$. Passing to a subsequence, we may assume that $(u_n)_{n \in \mathbb N}$ is convergent to a point $u_\ast$ and that $\norm{u_n-u_\ast}\leq \frac{1}{n^2}$ for all $n \in \mathbb N$. Then, $x:=(u_1,u_2,u_3,\ldots)$ belongs to $bv_1(E)$, but $C_f(x)\notin l^\infty(E)$. Therefore, $f$ must be bounded on any compact subset of $E$. 
\end{proof}

The completeness of the normed spaces $E$ in the above result turns out to be essential, as the following example shows.

\begin{example}\label{ex:bv_1_into_l_infty}
Set $a:=(1,\frac{1}{2^2}, \frac{1}{3^2},\ldots)$. Consider the map $f \colon c_{00}(\mathbb R) \to c_{00}(\mathbb R)$ given by the formula $f(u)=\frac{u}{\norm{u-a}_\infty}$. In Example~\ref{ex:lip_loc_vs_lip_comp} we saw that $f$ is Lipschitz continuous (and hence, bounded) on compact subsets of $c_{00}(\mathbb R)$. We also know that the sequence $x \colon N \to c_{00}(\mathbb R)$, given by $x(n):=P_n(a)$, belongs to $bv_1(c_{00}(\mathbb R))$ -- see Example~\ref{ex:completeness_essential}. But, $\norm{f(x(n))}_\infty=(n+1)^2$, which means that $C_f(x)\notin l^\infty(c_{00}(\mathbb R))$.
\end{example}

At this point it should  not be surprising that when $p>1$ we get completely different acting conditions for $C_f \colon bv_p(E) \to l^\infty(E)$.

\begin{theorem}\label{thm:bv_p_into_linfty}
Let $p \in (1,+\infty)$. Then, the following conditions are equivalent\textup:
\begin{enumerate}[label=\textup{(\roman*)}]
 \item\label{thm:bv_p_into_linfty_i} the composition operator $C_f$ maps $bv_p(E)$ into $l^\infty(E)$,
 \item\label{thm:bv_p_into_linfty_ii} $f$ is bounded on $E$,
 \item\label{thm:bv_p_into_linfty_iii} $f$ is bounded on countable subsets of $E$.
\end{enumerate}
\end{theorem}

\begin{proof}
As the implications $\ref{thm:bv_p_into_linfty_ii} \Rightarrow \ref{thm:bv_p_into_linfty_iii}$ and $\ref{thm:bv_p_into_linfty_iii} \Rightarrow \ref{thm:bv_p_into_linfty_i}$ are obvious, we need to prove  the implication $\ref{thm:bv_p_into_linfty_i} \Rightarrow \ref{thm:bv_p_into_linfty_ii}$ only. Suppose that $f$ is not bounded on $E$. Then, there exists a sequence $(u_n)_{n \in \mathbb N}$ of elements of $E$ such that $\sup_{n \in \mathbb N}\norm{f(u_n)}=+\infty$. Using the procedure described in Remark~\ref{rem:producing_sequences_in_bv_p}, we can construct $x \in bv_p(E)$ containing $(u_n)_{n\in \mathbb N}$ as its subsequence. But then, $C_f(x)\notin l^{\infty}(E)$.
\end{proof}

Finally, we turn our attention to composition operators acting into $c(E)$. Once again, we will distinguish two cases: $p=1$ and $p>1$.

\begin{theorem}\label{thm:bv_1_into_c}
Let $E$ be a Banach space. Then, the following conditions are equivalent\textup:
\begin{enumerate}[label=\textup{(\roman*)}]
 \item\label{thm:bv_1_into_c_i} the composition operator $C_f$ maps $bv_1(E)$ into $c(E)$,
 \item\label{thm:bv_1_into_c_ii} $f$ is continuous on $E$,
 \item\label{thm:bv_1_into_c_iii} $f$ is continuous on bounded subsets of $E$,
 \item\label{thm:bv_1_into_c_iv} $f$ is continuous on compact subsets of $E$,
 \item\label{thm:bv_1_into_c_v} $f$ is continuous on countable subsets of $E$.
\end{enumerate}
\end{theorem}

\begin{proof}
We will show that each of the first three conditions implies the next one, and the fourth implies the first again. With the last fifth condition we will deal separately.

$\ref{thm:bv_1_into_c_i} \Rightarrow \ref{thm:bv_1_into_c_ii}$ Suppose that there is a sequence $(u_n)_{n \in \mathbb N}$ of elements of $E$ which converges to a point $u_\ast \in E$ and for some $\varepsilon>0$ satisfies the inequality $\norm{f(u_n) - f(u_\ast)}\geq \varepsilon$ for $n \in \mathbb N$. Passing to a subsequence if necessary, we may assume that $\norm{u_n-u_\ast}\leq \frac{1}{n^2}$ for $n \in \mathbb N$. Now, set $x:=(u_\ast,u_1,u_\ast,u_2,\ldots)$. Clearly, $x \in bv_1(E)$. On the other hand, $\lim_{n \to \infty} f(x(2n-1))=f(u_\ast)$, whereas $\norm{f(x(2n))-f(u_\ast)}=\norm{f(u_n)-f(u_\ast)}\geq \varepsilon$ for $n \in \mathbb N$. This means that $C_f(x)\notin c(E)$.

The implications $\ref{thm:bv_1_into_c_ii} \Rightarrow \ref{thm:bv_1_into_c_iii}$ and $\ref{thm:bv_1_into_c_iii} \Rightarrow \ref{thm:bv_1_into_c_iv}$ are obvious. So, now we can turn to the implication $\ref{thm:bv_1_into_c_iv} \Rightarrow \ref{thm:bv_1_into_c_i}$. Assume that $f$ is continuous on compact subsets of $E$. Let $x \in bv_1(E)$. Then, $x \in c(E)$ and $K:=\{x_\ast\}\cup\dset{x(n) \in E}{n \in \mathbb N}$ is a compact subset of $E$; here $x_\ast$ is the limit of $x$. Consequently, $f(x(n)) \to f(x_\ast)$, by the assumption. Therefore, $C_f(x) \in c(E)$.

At this point we know that the first four conditions are equivalent. Hence, $\ref{thm:bv_1_into_c_i}$ implies $\ref{thm:bv_1_into_c_v}$. The proof of the other implication is straightforward, as it suffices to note that the set $K$ defined above is countable. 
\end{proof}

\begin{remark}
The fact that completeness of the normed space $E$ in Theorem~\ref{thm:bv_1_into_c} is essential follows from Example~\ref{ex:bv_1_into_l_infty} and the simple inclusion $c(E)\subset l^\infty(E)$. 
\end{remark}

What may come as a surprise is the fact that when $p>1$ (Lipschitz) continuous mappings, in general, do not generate composition operators that act from $bv_p(E)$ into $c(E)$. To see this it suffices to consider the identity mapping on $E$, that is, $f \colon E \to E$ given by $f(u)=u$. It generates the identity composition operator $C_f$ that maps $bv_p(E)$ onto $bv_p(E)$. But $bv_p(E) \not\subseteq c(E)$. And so, $C_f(bv_p(E))\not\subseteq c(E)$. In this case we need a significantly smaller class of generators.

\begin{theorem}\label{thm:bv_p_into_c}
Let $p \in (1,+\infty)$. The composition operator $C_f$ maps $bv_p(E)$ into $c(E)$ if and only if $f$ is a constant mapping.
\end{theorem}

\begin{proof}
Obviously, if $f$ is a constant mapping, then $C_f(bv_p(E))\subseteq c(E)$. So, now let us assume that there are two distinct points $u,w \in E$ such that $f(u)\neq f(w)$. Set $u_{2n-1}:=u$ and $u_{2n}:=w$ for $n \in \mathbb N$. Using the construction described in Remark~\ref{rem:producing_sequences_in_bv_p}, we can find $x \in bv_p(E)$ containing $(u_n)_{n\in \mathbb N}$ as its subsequence. But then, $C_f(x)$ contains two constant subsequences: one consisting of $f(u)$, and the other one consisting of $f(w)$. Thus, $C_f(x) \notin c(E)$. This ends the proof.
\end{proof}

\subsection{Acting conditions: summary}
\label{sec:summary_AC}
In this short section we summarize our results in a table. Let $E$ be a normed space. If $E$ is complete, we will highlight this fact by writing $\hat E$ instead of $E$. Further, given two normed spaces $X$ and $Y$ let $A(X,Y)$ be the set of all the maps $f \colon E \to E$ generating composition operators $C_f$ acting between $X$ and $Y$, that is, $A(X,Y):=\dset{f\colon E \to E}{\text{$C_f$ maps $X$ into $Y$}}$.

\begin{small}
\begin{table}[h!]
\begin{center}
\begin{tabular}{@{}ccc@{}}
\toprule
\parbox[c]{4cm}{\center $f \in A(c_0(E),bv_p(E))$\\ for $1 \leq p < +\infty$} & $\Leftrightarrow$ & there exists $\delta>0$ such that $f|_{\ball_E(0,\delta)}$ is constant \\[3mm] \midrule 
\parbox[c]{4cm}{\center $f \in A(c(E),bv_p(E))$\\ for $1\leq p < +\infty$} & $\Leftrightarrow$ & $f$ is constant\\[3mm] \midrule 
\parbox[c]{4cm}{\center $f \in A(l^\infty(E),bv_p(E))$\\ for $1\leq p < +\infty$} & $\Leftrightarrow$ & $f$ is constant\\[3mm] \midrule 
\parbox[c]{4cm}{\center $f \in A(l^p(E),bv_q(E))$\\ for $1\leq p,q < +\infty$} & $\Leftrightarrow$ & \parbox[c]{10.5cm}{\center there exist $\delta>0$ and $L\geq 0$ such that $\norm{f(u)-f(w)}\leq L\norm{u}^{\frac{p}{q}}+L\norm{w}^{\frac{p}{q}}$ for all $u,w \in \ball_E(0,\delta)$} \\[3mm] \midrule 
\parbox[c]{4cm}{\center $f \in A(bv_p(E),bv_q(E))$\\ for $1\leq q < p <+\infty$} & $\Leftrightarrow$ & $f$ is constant\\[3mm] \midrule 
\parbox[c]{4cm}{\center $f \in A(bv_1(\hat E),bv_p(\hat E))$\\ for $1\leq p <+\infty$} & $\Leftrightarrow$ & \parbox[c]{10.5cm}{\center for every compact subset $K$ of $\hat E$ there exists $L_K\geq 0$ such that $\norm{f(u)-f(w)}\leq L_K\norm{u-w}^{\frac{1}{p}}$ for all $u,w \in K$} \\[4mm] \midrule 
\parbox[c]{4cm}{\center $f \in A(bv_p(E),bv_q(E))$\\ for $1< p \leq q<+\infty$} & $\Leftrightarrow$ & \parbox[c]{10.5cm}{\center there exist $\delta>0$ and $L\geq 0$ such that $\norm{f(u)-f(w)}\leq L\norm{u-w}^{\frac{p}{q}}$ for all $u,w \in E$ with $\norm{u-w}\leq \delta$} \\[3.5mm] \midrule 
\parbox[c]{4cm}{\center $f \in A(bv_p(E),c_0(E))$\\ for $1 \leq p <+\infty$} & $\Leftrightarrow$ & $f$ is the zero map \\[3mm] \midrule 
\parbox[c]{4cm}{\center $f \in A(bv_p(E),l^q(E))$\\ for $1 \leq p,q <+\infty$} & $\Leftrightarrow$ & $f$ is the zero map \\[3mm] \midrule 
$f \in A(bv_1(\hat E),l^\infty(\hat E))$ & $\Leftrightarrow$ & $f$ is bounded on compact subsets of $\hat E$ \\ \midrule
\parbox[c]{4cm}{\center $f \in A(bv_p(E),l^\infty(E))$\\ for $1<p<+\infty$} & $\Leftrightarrow$ & $f$ is bounded on (countable subsets of) $E$ \\[3mm] \midrule
$f \in A(bv_1(\hat E),c(\hat E))$ & $\Leftrightarrow$ & $f$ is continuous on (bounded/compact/countable subsets of) $\hat E$\hspace{3mm} \\ \midrule
\parbox[c]{4cm}{\center $f \in A(bv_p(E),c(E))$\\ for $1<p<+\infty$} & $\Leftrightarrow$ & $f$ is constant\\[2mm]
\bottomrule
\end{tabular}
\vspace{1.5mm}
\caption{Acting conditions for the composition operator $C_f$.}
\end{center}
\end{table}
\end{small}

\clearpage

\section{Boundedness}
\label{sec:boundedness}

The main goal of this section is to characterize bounded and locally bounded composition operators acting to and from the space $bv_p(E)$. Let us recall that we will call a (non-linear) operator $F$ acting between two normed spaces $X$ and $Y$ \emph{bounded} if $F(X)$ is a bounded subset of $Y$. Similarly, we will call $F$ \emph{locally bounded},\label{page:locally_bounded} if $F(B)$ is a bounded subset of $Y$ for every bounded subset $B$ of $X$. Of course, when checking whether a given operator is locally bounded, it is enough to restrict all the considerations to closed balls of $X$. For completeness, let us mention that in the literature it is customary to include an additional assumption of continuity in the definition of a (locally) bounded operator. However, to be consisted with the conventions we adopted, and since we are going to discuss continuity of composition operators in $bv_p$-spaces in detail in the second part of our study, we decided to go for the ``non-continuous'' version of the definition. 

\subsection{Local boundedness}

Let us start with local boundedness. We will try to follow the same order of spaces as in Section~\ref{sec:action_conditions}. Certainly, we will exclude some of them from our discussion, as the composition operator acting between those spaces is generated by either a constant or the zero mapping.

\begin{theorem}\label{thm:bd_c0_bvp}
Let $p \in [1,+\infty)$ and assume that the composition operator $C_f$ maps $c_0(E)$ into $bv_p(E)$. Then, the following conditions are equivalent\textup:
\begin{enumerate}[label=\textup{(\roman*)}]
 \item\label{thm:bd_c0_bvp_i}  $C_f$ is bounded,
 \item\label{thm:bd_c0_bvp_ii} $C_f$ is locally bounded,
 \item\label{thm:bd_c0_bvp_iii} $f$ is a constant map.
\end{enumerate}
\end{theorem}

\begin{proof}
Clearly, only the implication~$\ref{thm:bd_c0_bvp_ii} \Rightarrow \ref{thm:bd_c0_bvp_iii}$ requires a proof. So, let us suppose that there is a point $u \in E$ such that $f(u)\neq f(0)$. Fix a positive integer $n$ and define $x_n \colon \mathbb N \to E$ as $(u,0,u,0,u,0,u,\ldots,u,0,0,0,\ldots)$ with the last non-zero term appearing on the $(2n-1)$-th position. Obviously, $x_n \in c_0(E)$ and $\norm{x_n}_\infty = \norm{u}$ for $n \in \mathbb N$. But
\begin{align*}
 \norm{C_f(x_n)}_{bv_p} \geq \Biggl(\sum_{k=1}^n \norm{f(u)-f(0)}^p \Biggr)^{\frac{1}{p}}=n^{\frac{1}{p}}\norm{f(u)-f(0)}.
\end{align*}   
This implies that $C_f \colon c_0(E) \to bv_p(E)$ is not locally bounded.
\end{proof}

Now, let us discuss local boundedness of composition operators acting from $l^p(E)$ into $bv_q(E)$. This time, we will not end up with constant generators. 

\begin{theorem}\label{thm:local_boundedness_lp_bvq}
Let $p,q \in [1,+\infty)$ and assume that the composition operator $C_f$ maps $l^p(E)$ into $bv_q(E)$. Then, $C_f$ is locally bounded if and only if $f$ is.
\end{theorem}

\begin{proof}
First, we will prove the necessity part, i.e., we will show that $f$ is locally bounded provided that $C_f$ is. To this end let us fix $r>0$ and choose $R>0$ so that $C_f(\ball_{l^p}(0,r))\subseteq \ball_{bv_q}(0,R)$. Now, take $u \in \ball_E(0,r)$. Also, define $x:=(u,0,0,0,\ldots)$. Then, $x \in \ball_{l^p}(0,r)$. And hence, $\norm{f(u)}=\norm{f(x(1))} = \norm{C_f(x)(1)}\leq \norm{C_f(x)}_{bv_q}\leq R$. This shows that $f$ is locally bounded.

Now, let us move to the second part of the proof. This time, we assume that $f$ is locally bounded. Let us fix $r>0$. Then, there is $\rho>0$ such that $f(\ball_E(0,r))\subseteq \ball_E(0,\rho)$. Further, since $C_f$ maps $l^p(E)$ into $bv_q(E)$, by Theorem~\ref{thm:lp_into_bvq} there exist $\delta>0$ and $L\geq 0$ such that $\norm{f(u)-f(w)}\leq L\norm{u}^{\frac{p}{q}}+L\norm{w}^{\frac{p}{q}}$ for all $u,w \in \ball_E(0,\delta)$. Set
\[
 R:=\rho + 2\bigl(2L^q r^p + 4r^p \rho^q \delta^{-p}\bigr)^{\frac{1}{q}}.
\]   
Take $x \in \ball_{l^p}(0,r)$ and consider the following (pairwise) disjoint sets of indices
\begin{align*}
 I_1&:=\dset[\big]{n \in \mathbb N}{\text{$\norm{x(n)}\leq \delta$ and $\norm{x(n+1)}\leq \delta$}},\\
  I_2&:=\dset[\big]{n \in \mathbb N}{\text{$\norm{x(n)}> \delta$ and $\norm{x(n+1)}\leq \delta$}},\\
 I_3&:=\dset[\big]{n \in \mathbb N}{\text{$\norm{x(n)}\leq \delta$ and $\norm{x(n+1)}> \delta$}},\\
 I_4&:=\dset[\big]{n \in \mathbb N}{\text{$\norm{x(n)}> \delta$ and $\norm{x(n+1)}> \delta$}}.
\end{align*}
Note that only the set $I_1$ is infinite. If by $\abs{J}$ we denote the number of elements in $J$, then
\begin{align*}
2r^p \geq 2 \sum_{n=1}^\infty \norm{x(n)}^p \geq \sum_{j=2}^4 \sum_{n \in I_j} \bigl(\norm{x(n+1)}^p + \norm{x(n)}^p\bigr)\geq  \sum_{j=2}^4 \abs{I_j}\cdot \delta^p =\abs{\mathbb N\setminus I_1} \cdot \delta^p.
\end{align*}
This, in turn, implies that
\begin{align*}
&\sum_{n=1}^\infty \norm[\big]{f(x(n+1))-f(x(n))}^q\\
&\qquad = \sum_{n \in I_1} \norm[\big]{f(x(n+1))-f(x(n))}^q + \sum_{n \in \mathbb N\setminus I_1} \norm[\big]{f(x(n+1))-f(x(n))}^q\\
&\qquad \leq L^q \cdot \sum_{n \in I_1}\bigl( \norm{x(n+1)}^{\frac{p}{q}}+\norm{x(n)}^{\frac{p}{q}}\bigr)^q + \abs{\mathbb N \setminus I_1}\cdot 2^{q+1}\rho^q\\
&\qquad \leq 2^{q+1}L^q\cdot \sum_{n=1}^\infty \norm{x(n)}^p + 2^{q+2} r^p \rho^q \delta^{-p}. 
\end{align*}
Thus, $\norm{C_f(x)}_{bv_q}\leq R$. This ends the proof. 
\end{proof}

For the composition operator acting between $bv_p(E)$ and $bv_q(E)$, where $1<p\leq q<+\infty$, local boundedness is a direct consequence of the acting conditions.

\begin{theorem}\label{thm:AC_implies_B_loc_bvp}
Let $1<p \leq q<+\infty$. Moreover, assume that the composition operator $C_f$ maps $bv_p(E)$ into $bv_q(E)$. Then, $C_f$ is locally bounded. 
\end{theorem}

\begin{proof}
Fix $r>0$. Since $C_f$ maps $bv_p(E)$ into $bv_q(E)$, by Proposition~\ref{prop:eq_lip_strg} there exists $L_r\geq 0$ such that $\norm{f(u)-f(w)}\leq L_r\norm{u-w}^{\frac{p}{q}}$ for all $u,w \in E$ with $\norm{u-w}\leq r$. Set $R:=\norm{f(0)}+2L_r r^{\frac{p}{q}}$. Now, take any $x \in \ball_{bv_p}(0,r)$. Then, $\norm{x(1)}\leq r$ and $\norm{x(n+1)-x(n)}\leq r$ for $n \in \mathbb N$. This, in turn, implies that
\[
 \norm{C_f(x)(1)} \leq \norm{f(0)}+\norm{f(x(1))-f(0)} \leq \norm{f(0)} + L_r\norm{x(1)}^{\frac{p}{q}}\leq \norm{f(0)}+L_r r^{\frac{p}{q}}
\]
and
\[
 \norm{C_f(x)(n+1)-C_f(x)(n)}^q = \norm{f(x(n+1))-f(x(n))}^q \leq L_r^q \norm{x(n+1)-x(n)}^p
\]
for $n \in \mathbb N$. Thus, $\norm{C_f(x)}_{bv_q} \leq R$, which means that $C_f$ is locally bounded.
\end{proof}

Looking at the previous theorem and the corresponding results in the parallel theory of composition operators in the spaces $BV_p([a,b],\mathbb R)$ (cf.~\cite{ABM}*{Chapter~5} or~Section~\ref{sec:comparision_with_BV} below), it would be quite natural to expect that also for composition operators acting between $bv_1(E)$ and $bv_p(E)$ acting conditions imply local boundedness. However, this is not the case. In the following example we will define a bounded map $f$ that is Lipschitz continuous on compact sets, but does not generate a locally bounded composition operator. 

\begin{example}
Let us consider the bounded map $f \colon l^2(\mathbb R) \to l^2(\mathbb R)$ defined in Example~\ref{ex:lip_loc_vs_lip_comp_2}. We know that $f$ is Lipschitz continuous on compact subsets of $l^2(\mathbb R)$, and so it generates the composition operator $C_f$ that acts from $bv_1(l^2(\mathbb R))$ into itself. However, as we are going to show, $C_f$ is not locally bounded. By $e_k$ let us denote the $k$-th unit vector of $l^2(\mathbb R)$. For each $n \in \mathbb N$ set 
\[
x_n:=\bigl(e_n, (1-n^{-2})e_n, e_n, (1-n^{-2})e_n, \ldots, e_n, (1-n^{-2})e_n,0,0,0,\ldots\bigr),
\]
where the last non-zero element appears at the $2n^2$-th position. Then, $x_n \in bv_1(l^2(\mathbb R))$ and $\norm{x_n}_{bv_1} \leq 1+(2n^2-1)\norm{e_n -(1-n^{-2})e_n}_{l^2} + 1-n^{-2} \leq 4$ for $n \in \mathbb N$. Since $f(e_n)=0$ and $f\bigl((1-n^{-2})e_n\bigr)=(-\sin[(2\pi+1)n^{-1}],0,0,\ldots)$, for $n \geq 5$ we thus have
\[
\norm{C_f(x_n)}_{bv_1} \geq n^2 \norm[\big]{f(e_n) - f\bigl((1-\tfrac{1}{n^{2}})e_n\bigr)}_{l^2} \geq 4n
\]
(cf.~Example~\ref{ex:lip_loc_vs_lip_comp_2}). This means that $C_f$ is not locally bounded.
\end{example}

There is a very simple reason why the composition operator $C_f$ in the above example fails to be locally bounded. Its generator $f$ is not Lipschitz continuous on bounded subsets of $l^2(\mathbb R)$ -- see Example~\ref{ex:lip_loc_vs_lip_comp_2}. Therefore, it should not come as a surprise that H\"older continuity on bounded sets of $f$ is equivalent to local boundedness of $C_f \colon bv_1(E) \to bv_p(E)$.

\begin{theorem}\label{thm:B_loc_bv1_bvp}
Let $p \in [1,+\infty)$ and let $E$ be a Banach space. Moreover, assume that the composition operator $C_f$ maps $bv_1(E)$ into $bv_p(E)$. Then, $C_f$ is locally bounded if and only if $f$ is H\"older continuous on bounded subsets of $E$ with exponent $\frac{1}{p}$.
\end{theorem}

\begin{proof}
First, let us assume that $f$ is H\"older continuous on bounded subsets of $E$ with exponent $\frac{1}{p}$. Also, fix $r>0$. Then, there exists a constant $L_r\geq 0$ such that $\norm{f(u)-f(w)}\leq L_r\norm{u-w}^{\frac{1}{p}}$ for any $u,w \in \ball_E(0,r)$. Set $R:=2L_r r^{\frac{1}{p}} + \norm{f(0)}$. If $x \in \ball_{bv_1}(0,r)$, then $x(n) \in \ball_E(0,r)$ for $n \in \mathbb N$. Hence, $\norm{f(x(n+1))-f(x(n))}^p \leq L_r^p\norm{x(n+1)-x(n)}$ for all $n \in \mathbb N$. Furthermore, $\norm{f(x(1))}\leq \norm{f(x(1))-f(0)}+\norm{f(0)} \leq L_r\norm{x(1)}^{\frac{1}{p}}+\norm{f(0)}$. All of this implies that
\begin{align*}
 \norm{C_f(x)}_{bv_p} &= \norm{f(x(1))} + \Biggl(\sum_{n=1}^\infty \norm{f(x(n+1))-f(x(n))}^p \Biggr)^{\frac{1}{p}}\\
& \leq L_r\norm{x(1)}^{\frac{1}{p}}+\norm{f(0)} + \Biggl(\sum_{n=1}^\infty L_r^p\norm{x(n+1)-x(n)} \Biggr)^{\frac{1}{p}}\\
& \leq  L_r r^{\frac{1}{p}} + \norm{f(0)} + L_r r^{\frac{1}{p}}=R.
\end{align*}
Therefore, $C_f$ is locally bounded.

Now, let us proceed to the second part of the proof. Let us underline that some of the estimates will be only sketched, as similar ones were given in full detail in the proof of Proposition~\ref{prop:josephy_bv}. Assume that $C_f$ maps $bv_1(E)$ into $bv_p(E)$ and is locally bounded. This, in particular, implies that $f$ is locally bounded as well (cf. the proof of Theorem~\ref{thm:local_boundedness_lp_bvq}). So, $M_\rho:=\sup_{\norm{u}\leq \rho} \norm{f(u)}<+\infty$ for every $\rho>0$. Suppose now that $f$ is not H\"older continuous with exponent $\frac{1}{p}$ on some ball $\ball_E(0,r)$. Then, there exist sequences $(u_n)_{n \in \mathbb N}$ and $(w_n)_{n \in \mathbb N}$ of elements of $\ball_E(0,r)$ such that
\[
 \norm{f(u_n)-f(w_n)} > 2^{1+\frac{1}{p}}M_r n^3 \norm{u_n-w_n}^{\frac{1}{p}}.
\]
For each $n \in \mathbb N$ set $m_n:=[n^{-2p}\norm{u_n-w_n}^{-1}]$ and $x_n:=(u_n,w_n,u_n,w_n,\ldots,u_n,w_n,0,0,\ldots)$, where the last element $w_n$ appears at the $2m_n$-th position and $[t]$ denotes the greatest integer smaller than or equal to $t$. Then, clearly, $x_n \in bv_1(E)$ and
\[
 \norm{x_n}_{bv_1} = \norm{u_n} + (2m_n-1)\norm{u_n-w_n} + \norm{w_n} \leq 2r + 2n^{-2p} \leq 2+2r
\]
for $n \in \mathbb N$. However, 
\[
\text{$\norm{C_f(x_n)}^p_{bv_p} \geq m_n\norm{f(u_n)-f(w_n)}^p \geq (2M_rn)^p \to +\infty$ as $n \to +\infty$.}
\]
This shows that $C_f$ cannot be locally bounded -- a contradiction.
\end{proof}

We end this subsection with a result on local boundedness of composition operators acting between $bv_1(E)$ and $c(E)$ or $l^\infty(E)$. We will, however, omit its proof, as it follows directly from the definition, and should be evident at this point. 

\begin{theorem}\label{thm:local_bounded_bv1_to_c}
Let $E$ be a Banach space and let $X \in \{c(E),l^\infty(E)\}$. Moreover, assume that the composition operator $C_f$ maps $bv_1(E)$ into $X$. Then, $C_f$ is locally bounded if and only if $f$ is. 
\end{theorem}

\subsection{Boundedness}

In this very short section we focus on the boundedness of composition operators acting in the sequence spaces of our interest. Without further ado let us state the first result. 

\begin{theorem}\label{thm:boundedness}
Let $p,q \in [1,+\infty)$ and let $X\in \{l^p(E),bv_p(E)\}$. Moreover, assume that the composition operator $C_f$ maps $X$ into $bv_q(E)$. Then, $C_f$ is bounded if and only if $f$ is a constant map.
\end{theorem}

\begin{proof}
Clearly, we need to prove only the necessity part. So, suppose that $f$ is not a constant map, that is, there exists $u \in E$ such that $f(u)\neq f(0)$. For each $n \in \mathbb N$ define $x_n:=(u,0,u,0,\ldots,u,0,0,\ldots)$, where the last non-zero term appears on the $(2n-1)$-th position. Then,
\[
\norm{C_f(x)}_{bv_q} \geq  \Biggl(\sum_{k=1}^n \norm{f(u)-f(0)}^q \Biggr)^{\frac{1}{q}}=n^{\frac{1}{q}}\norm{f(u)-f(0)}
\]
for $n \in \mathbb N$. This shows that $C_f$ is not bounded.  
\end{proof}

On a more positive note, there are non-trivial bounded composition operators acting between $bv_1(E)$ and $c(E)$ or $l^\infty(E)$.

\begin{theorem}
Let $E$ be a Banach space and let $X \in \{c(E),l^\infty(E)\}$. Moreover, assume that the composition operator $C_f$ maps $bv_1(E)$ into $X$. Then, $C_f$ is bounded if and only if $f$ is. 
\end{theorem}

\begin{proof}
Clearly, we need to prove only the necessity part. Suppose that $f$ is not bounded on $E$. Then, there is a sequence $(u_n)_{n \in \mathbb N}$ of elements of $E$ such that $\sup_{n \in \mathbb N}\norm{f(u_n)}=+\infty$. For each $n \in \mathbb N$ set $x_n:=(u_n,0,0,\ldots)$. Then, $x_n \in bv_1(E)$ and $\norm{C_f(x_n)}_\infty\geq \norm{f(x_n(1))}=\norm{f(u_n)}$ for $n \in \mathbb N$. Therefore, $C_f$ is not bounded. 
\end{proof}

\subsection{Boundedness: summary}

Similar to what we did in Section~\ref{sec:summary_AC}, here, we will summarize conditions on the generator $f$ that guarantee the (local) boundedness of $C_f$. Given two normed spaces $X$ and $Y$, by $B_{\text{loc}}(X,Y)$ and $B(X,Y)$ we will denote the sets of all the maps $f \colon E \to E$ generating composition operators $C_f \colon X \to Y$ that are locally bounded and bounded, respectively. Let us recall that if the normed space $E$ is complete, we will underline this fact by writing $\hat E$ instead of $E$.

\begin{small}
\begin{table}[hbt!]
\begin{center}
\begin{tabular}{@{}ccc@{}}
\toprule
\parbox[c]{4.2cm}{\center $f \in B_{\text{loc}}(c_0(E),bv_p(E))$\\ for $1 \leq p < +\infty$} & $\Leftrightarrow$ & $f$ is constant\\[3mm] \midrule 
\parbox[c]{4.2cm}{\center $f \in B_{\text{loc}}(c(E),bv_p(E))$\\ for $1\leq p < +\infty$} & $\Leftrightarrow$ & $f$ is constant\\[3mm] \midrule 
\parbox[c]{4.2cm}{\center $f \in B_{\text{loc}}(l^\infty(E),bv_p(E))$\\ for $1\leq p < +\infty$} & $\Leftrightarrow$ & $f$ is constant\\[3mm] \midrule 
\parbox[c]{4.2cm}{\center $f \in B_{\text{loc}}(l^p(E),bv_q(E))$\\ for $1\leq p,q < +\infty$} & $\Leftrightarrow$ & \parbox[c]{10.5cm}{\center $f$ is locally bounded and there exist $\delta>0$ and $L\geq 0$ such that $\norm{f(u)-f(w)}\leq L\norm{u}^{\frac{p}{q}}+L\norm{w}^{\frac{p}{q}}$ for all $u,w \in \ball_E(0,\delta)$} \\[3mm] \midrule 
\parbox[c]{4.2cm}{\center $f \in B_{\text{loc}}(bv_p(E),bv_q(E))$\\ for $1\leq q < p <+\infty$} & $\Leftrightarrow$ & $f$ is constant\\[3mm] \midrule 
\parbox[c]{4.2cm}{\center $f \in B_{\text{loc}}(bv_1(\hat E),bv_p(\hat E))$\\ for $1\leq p <+\infty$} & $\Leftrightarrow$ & \parbox[c]{10.5cm}{\center for every $r>0$ there exists $L_r\geq 0$ such that $\norm{f(u)-f(w)}\leq L_r\norm{u-w}^{\frac{1}{p}}$ for all $u,w \in \ball_{\hat{E}}(0,r)$} \\[4mm] \midrule 
\parbox[c]{4.2cm}{\center $f \in B_{\text{loc}}(bv_p(E),bv_q(E))$\\ for $1< p \leq q<+\infty$} & $\Leftrightarrow$ & \parbox[c]{10.5cm}{\center there exist $\delta>0$ and $L\geq 0$ such that $\norm{f(u)-f(w)}\leq L\norm{u-w}^{\frac{p}{q}}$ for all $u,w \in E$ with $\norm{u-w}\leq \delta$} \\[3.5mm] \midrule 
\parbox[c]{4.2cm}{\center $f \in B_{\text{loc}}(bv_p(E),c_0(E))$\\ for $1 \leq p <+\infty$} & $\Leftrightarrow$ & $f$ is the zero map \\[3mm] \midrule 
\parbox[c]{4.2cm}{\center $f \in B_{\text{loc}}(bv_p(E),l^q(E))$\\ for $1 \leq p,q <+\infty$} & $\Leftrightarrow$ & $f$ is the zero map \\[3mm] \midrule 
$f \in B_{\text{loc}}(bv_1(\hat E),l^\infty(\hat E))$ & $\Leftrightarrow$ & $f$ is locally bounded\\ \midrule
\parbox[c]{4.2cm}{\center $f \in B_{\text{loc}}(bv_p(E),l^\infty(E))$\\ for $1<p<+\infty$} & $\Leftrightarrow$ & $f$ is bounded \\[3mm] \midrule
$f \in B_{\text{loc}}(bv_1(\hat E),c(\hat E))$ & $\Leftrightarrow$ & $f$ is locally bounded and continuous on $\hat E$\hspace{3mm} \\ \midrule
\parbox[c]{4.2cm}{\center $f \in B_{\text{loc}}(bv_p(E),c(E))$\\ for $1<p<+\infty$} & $\Leftrightarrow$ & $f$ is constant \\[2mm]
\bottomrule
\end{tabular}
\vspace{1.5mm}
\caption{Local boundedness of the composition operator $C_f$.}
\end{center}
\end{table}
\end{small}

\clearpage

\begin{small}
\begin{table}[hbt!]
\begin{center}
\begin{tabular}{@{}ccc@{}}
\toprule
\parbox[c]{4.2cm}{\center $f \in B(c_0(E),bv_p(E))$\\ for $1 \leq p < +\infty$} & $\Leftrightarrow$ & $f$ is constant\\[3mm] \midrule 
\parbox[c]{4.2cm}{\center $f \in B(c(E),bv_p(E))$\\ for $1\leq p < +\infty$} & $\Leftrightarrow$ & $f$ is constant\\[3mm] \midrule 
\parbox[c]{4.2cm}{\center $f \in B(l^\infty(E),bv_p(E))$\\ for $1\leq p < +\infty$} & $\Leftrightarrow$ & $f$ is constant\\[3mm] \midrule 
\parbox[c]{4.2cm}{\center $f \in B(l^p(E),bv_q(E))$\\ for $1\leq p,q < +\infty$} & $\Leftrightarrow$ & $f$ is constant\\[3mm] \midrule 
\parbox[c]{4.2cm}{\center $f \in B(bv_p(E),bv_q(E))$\\ for $1\leq q < p <+\infty$} & $\Leftrightarrow$ & $f$ is constant\\[3mm] \midrule 
\parbox[c]{4.2cm}{\center $f \in B(bv_1(\hat E),bv_p(\hat E))$\\ for $1\leq p <+\infty$} & $\Leftrightarrow$ & $f$ is constant\\[4mm] \midrule 
\parbox[c]{4.2cm}{\center $f \in B(bv_p(E),bv_q(E))$\\ for $1< p \leq q<+\infty$} & $\Leftrightarrow$ & $f$ is constant\\[3.5mm] \midrule 
\parbox[c]{4.2cm}{\center $f \in B(bv_p(E),c_0(E))$\\ for $1 \leq p <+\infty$} & $\Leftrightarrow$ & $f$ is the zero map \\[3mm] \midrule 
\parbox[c]{4.2cm}{\center $f \in B(bv_p(E),l^q(E))$\\ for $1 \leq p,q <+\infty$} & $\Leftrightarrow$ & $f$ is the zero map \\[3mm] \midrule 
$f \in B(bv_1(\hat E),l^\infty(\hat E))$ & $\Leftrightarrow$ & $f$ is bounded\\ \midrule
\parbox[c]{4.2cm}{\center $f \in B(bv_p(E),l^\infty(E))$\\ for $1<p<+\infty$} & $\Leftrightarrow$ & $f$ is bounded \\[3mm] \midrule
$f \in B(bv_1(\hat E),c(\hat E))$ & $\Leftrightarrow$ & $f$ is bounded and continuous on $\hat E$\hspace{3mm} \\ \midrule
\parbox[c]{4.2cm}{\center $f \in B(bv_p(E),c(E))$\\ for $1<p<+\infty$} & $\Leftrightarrow$ & $f$ is constant\\[2mm]
\bottomrule
\end{tabular}
\vspace{1.5mm}
\caption{Boundedness of the composition operator $C_f$.}
\end{center}
\end{table}
\end{small}

\newpage

\section{Discussion and conclusions}
\label{sec:comparison}

In this last section, we will compare the results we obtained with some other theorems on composition/superposition operators in sequences spaces that are already known in the literature. We will also discuss how our findings relate to and differ from those in the parallel theory of composition operators in $BV_p$-spaces.

\subsection{Composition operators in $bv_p$- and $BV_p$-spaces}
\label{sec:comparision_with_BV}

The theory of composition and superposition operators in spaces of functions of bounded variation is very rich (see, for example,~\cite{ABM}*{Chapters~5 and~6},~\cite{DN}*{Chapter~6} and~\cite{reinwand}*{Chapter~5}). At its center lies the Josephy theorem, which, for spaces of Wiener variation, asserts that a composition operator $C_f$ maps $BV_p([a,b],E)$ into $BV_q([a,b],E)$, where $1\leq p\leq q<+\infty$, if and only if its generator $f \colon E \to E$ is H\"older continuous on precompact subsets of $E$ with exponent $\frac{p}{q}$ (see~\cite{BK}*{Theorem~4}); here, as in previous sections, $E$ stands for the normed space over the field of either real or complex numbers. In complete metric spaces, precompact subsets are exactly those whose closures (with respect to the whole space) are compact. Therefore, Josephy's result matches up with our findings when $p=1$ and $q \in [1,+\infty)$ -- see Theorem~\ref{thm:1}. However, when $1<p\leq q<+\infty$, the requirements for a map $ f \colon E \to E$ to generate the composition operator $C_f$ that acts between $bv_p(E)$ and $bv_q(E)$ are much stronger compared to the $BV_p$-setting (see Section~\ref{sec:generators} and Theorem~\ref{thm:bv_p_int_bv_q_strong}). The reason for this is simple: $bv_p(E)$, where $p>1$, contains unbounded sequences, while all functions in $BV_p([a,b],E)$ are bounded.

One of the consequences of Josephy's theorem in the real-valued setting is the local boundedness principle. It states that all composition operators mapping $BV_p([a,b], \mathbb R)$ into $BV_q([a,b], \mathbb R)$ are locally bounded (see, for example,~\cite{ABM}*{Theorem~5.27}). A similar result for sequences of bounded variation is true only for composition operators $C_f$ between $bv_p(E)$ and $bv_q(E)$ for $1<p\leq q < +\infty$; here, the normed space $E$ may be of either finite or infinite dimension -- see Theorem~\ref{thm:AC_implies_B_loc_bvp}. However, when we move to the abstract setting, the acting conditions for $C_f \colon BV_p([a,b],E) \to BV_q([a,b],E)$ no longer imply the local boundedness of the composition operator. It turns out that $C_f \colon BV_p([a,b],E) \to BV_q([a,b],E)$, where $1\leq p \leq q <+\infty$, is locally bounded if and only if its generator $f \colon E \to E$ is H\"older continuous on bounded subsets of $E$ with exponent $\frac{p}{q}$ (see~\cite{BK}*{Theorem~6} and cf.~\cite{DN}*{Proposition~6.34}). This result aligns with our findings for $C_f \colon bv_p(E) \to bv_q(E)$, when $p=1$ and $q \in [1,+\infty)$ -- see Theorem~\ref{thm:B_loc_bv1_bvp}.

When it comes to the boundedness of composition operators in $bv_p$- and $BV_p$-spaces, the situation is very simple. In both cases, the only bounded composition operators are the constant ones (see Theorem~\ref{thm:boundedness} and~\cite{BK}*{Remark~7}).

\subsection{Superposition operators in $bv_p$-spaces}
\label{sec:Karami}

The main motivation for our study comes from the recent paper~\cite{KFA}. While it deals with various properties of superposition operators in $bv_p$-spaces, here, we will focus on the acting conditions only. The reason behind our choice is simple: other results share similar flaws.

On page~311 of~\cite{KFA}, the following characterization of superposition operators acting between $bv_p(\mathbb R)$ and $bv_q(\mathbb R)$ is presented. Note that, for consistency, we adopt ``our'' notation. Additionally, it is important to mention that all the sequences considered in~\cite{KFA} start with the zero term. The normed subspace of $bv_p(\mathbb R)$ consisting of all such sequences will be denoted by $bv_p^0(\mathbb R)$.

\begin{theorema}\label{thmA}
Let $1\leq p,q < +\infty$ and let $S_f$ be the superposition operator generated by a function $f \colon \mathbb N \times \mathbb R \to \mathbb R$ such that $f(n,0)=0$ for $n \in \mathbb N$. Then, the following conditions are equivalent\textup:
\begin{enumerate}[label=\textup{(\alph*)}]
  \item\label{thmA_a} $S_f$ maps $bv_p^0(\mathbb R)$ into $bv_q^0(\mathbb R)$,

 \item\label{thmA_b} for every $x \in bv_p^0(\mathbb R)$ there exist a sequence $a \in bv_q^0(\mathbb R)$ and constants $\delta>0$, $L\geq 0$ and $N \in \mathbb N$ such that
\[
 \abs{f(n+1,x(n+1))-f(n,x(n))} \leq a(n+1)-a(n) + L\abs{x(n+1)-x(n)}^{\frac{p}{q}}
\]
for $n \geq N$ and $\abs{x(n+1)-x(n)}<\delta$,

 \item\label{thmA_c} for every $x \in bv_p^0(\mathbb R)$ and $\varepsilon>0$ there exist a sequence $a_\varepsilon \in bv_q^0(\mathbb R)$ and constants $\delta_\varepsilon>0$, $L_\varepsilon\geq 0$ and $N_\varepsilon \in \mathbb N$ such that
\[
 \abs{f(n+1,x(n+1))-f(n,x(n))} \leq a_\varepsilon(n+1)-a_\varepsilon(n) + L_\varepsilon\abs{x(n+1)-x(n)}^{\frac{p}{q}}
\]
for $n \geq N_\varepsilon$ and $\abs{x(n+1)-x(n)}<\delta_\varepsilon$.
\end{enumerate}
\end{theorema}

A couple of issues regarding Theorem~\ref{thmA} need to be raised. Firstly, contrary to what is usually expected, the conditions~\ref{thmA_b} and~\ref{thmA_c} are actually expressed in terms of the superposition operator $S_f$ rather than its generator, as they must hold for all sequences $x \in bv_p^0(\mathbb R)$. This might make their verification difficult in specific situations.  

Secondly, it might be somewhat surprising that a single result can cover various cases, regardless of the values of the parameters $p$ and $q$. In Section~\ref{sec:AC_into_bv}, we saw that depending on whether $1=p\leq q < +\infty$, $1<p\leq q < +\infty$, or $1\leq q < p < +\infty$, the classes of generators of the composition operators acting between the spaces $bv_p(E)$ and $bv_q(E)$ are essentially different.

Finally, let us add that certain steps in the proof of Theorem~\ref{thmA} are likely incorrect for $p>1$. Given $x \in bv_p^0(\mathbb R)$ such that $\norm{x}_{bv_p}\leq r$ for some $r>0$, we cannot generally claim that $\norm{(P_n-P_m)x}_{bv_p}\leq \varphi(r)$ without controlling the distance between the indices $n$ and $m$; here, $\varphi \colon (0,+\infty) \to (0,+\infty)$ is a given function of one variable. To see this, consider the sequence $x=(0,\frac{1}{k}, \frac{1}{k}+\frac{1}{k+1},\frac{1}{k}+\frac{1}{k+1}+\frac{1}{k+2},\ldots)$, where $k \in \mathbb N$ is chosen so that $\sum_{j=k}^\infty j^{-p} \leq r^p$. Also, fix $n \geq 2$ and observe that while $\norm{x}_{bv_p}\leq r$, we have $\norm{(P_n-P_m)x}_{bv_p} \geq \sum_{j=k}^{k+m-2} j^{-1} \to +\infty$ as $m \to +\infty$. 

\subsection{Superposition operators in $l^p$-spaces}

One of the classical papers on superposition operators in $l^p$-spaces was written by Dedagich and Zabre\u{\i}ko (see~\cite{DZ}). Since $bv_p(E)$ and $E\times l^p(E)$ are linearly isometric for $p\geq 1$, a natural question arises: can the ideas from~\cite{DZ} be applied to the $bv_p$-spaces? Well, the answer is both yes and no. While understanding the composition/superposition operators in $l^p$-spaces can help figure out the right classes of generators for their $bv_p$-counterparts, it seems that there is no straightforward way to translate all results from one setting to the other. It does not mean, however, that no connections can be found. To better understand what happens, let us consider a theoretical example.    

A corollary to Theorem~1 in~\cite{DZ} characterizes composition operators acting in $l^p$-spaces and reads as follows.

\begin{theorema}\label{thm:b}
Let $p,q \in [1,+\infty)$ and let $E$ be a normed space. A composition operator $C_f$, generated by a function $f\colon E \to E$, maps $l^p(E)$ into $l^q(E)$ if and only if
there exist constants $\delta>0$ and $L\geq 0$ such that $\norm{f(u)}\leq L\norm{u}^{\frac{p}{q}}$ for all $u \in \ball_E(0,\delta)$.
\end{theorema}

Clearly, Theorem~\ref{thm:b} may be viewed as an archetype for ``our'' Theorem~\ref{thm:lp_into_bvq}. On the other hand, it is Theorem~\ref{thm:lp_into_bvq}'s direct consequence. Thus, we obtain an alternative method for proving the characterization of acting conditions for $C_f \colon l^p(E)\to l^q(E)$, which is based on a completely different approach than the one used in~\cite{DZ}. As a final note, let us mention that the results obtained by Dedagich and Zabre\u{\i}ko were limited to the real-valued case; this can make a difference, especially if the classes of suitable generators revolve around H\"older continuous mappings -- for a detailed discussion of this phenomenon in the case of $BV_p$-spaces see~\cite{BK}.

And now, for completes, let us present the proof of Theorem~\ref{thm:b} based on Theorem~\ref{thm:lp_into_bvq}.

\begin{proof}[Proof of Theorem~\ref{thm:b}]
Clearly, only the necessity requires a proof. So, let us assume that $C_f$ maps $l^p(E)$ into $l^q(E)$. Since $l^q(E)\subseteq bv_q(E)$, in view of Theorem~\ref{thm:lp_into_bvq}, there exist two constants $\delta>0$ and $L\geq 0$ such that $\norm{f(u)-f(0)}\leq L\norm{u}^{\frac{p}{q}}$ for all $u \in \ball_E(0,\delta)$. To end the poof note that $f(0)=0$, because the only constant sequence in $l^q(E)$ is the zero one.
\end{proof}

It would be equally tempting to ``reverse'' the reasoning and apply Theorem~\ref{thm:b} to proving characterization of composition operators between $bv_p(E)$ and $bv_q(E)$.
It seems that this cannot be easily done for at least two reasons. Firstly, composition operators are non-linear mappings and generally do not preserve linear structure or isometries. And secondly, as we saw in Section~\ref{sec:AC_into_bv}, depending on the values of the parameters $p$ and $q$, there are essentially three different results describing acting conditions for $C_f \colon bv_p(E) \to bv_q(E)$.

\end{document}